%% file: Compositionality.tex
\documentclass[11pt]{article} 
\usepackage{amsthm}
\usepackage{amsmath,amssymb}       

\usepackage{xspace}
\usepackage{color}
\usepackage{epsfig}    
\graphicspath{{.}{./figures/}}  

\usepackage{macros/tikzfig}
\usepackage{keycommand}

\input{macros/defs.tex}

\input{macros/tikzstyles.tex}

\input{macros/tikzfigures.tex}
\definecolor{hexcolor0xa9a9a9}{rgb}{0.663,0.663,0.663}
\tikzstyle{GrayLine}=[dashed,draw=hexcolor0xa9a9a9] 
\tikzstyle{gray}=[dashed,draw=hexcolor0xa9a9a9]

\theoremstyle{definition}
\newtheorem{theorem}{Theorem}[section]
\newtheorem*{theorem*}{Theorem}

\newtheorem{defn}[theorem]{Definition}

\newtheorem{example*}[theorem]{Example*}
\newtheorem{examples*}[theorem]{Examples*}

\newtheorem{remark*}[theorem]{Remark*}

\def\bR{\begin{color}{red}}  
\def\bB{\begin{color}{blue}}
\def\bM{\begin{color}{magenta}}  
\def\bC{\begin{color}{cyan}}
\def\bW{\begin{color}{white}}
\def\bBl{\begin{color}{black}}
\def\bG{\begin{color}{green}}
\def\bY{\begin{color}{yellow}}
\def\e{\end{color}\xspace}
\newcommand{\bit}{\begin{itemize}}
\newcommand{\eit}{\end{itemize}\par\noindent}
\newcommand{\ben}{\begin{enumerate}}
\newcommand{\een}{\end{enumerate}\par\noindent}
\newcommand{\beq}{\begin{equation}}
\newcommand{\eeq}{\end{equation}\par\noindent}
\newcommand{\beqa}{\begin{eqnarray*}}
\newcommand{\eeqa}{\end{eqnarray*}\par\noindent}
\newcommand{\beqn}{\begin{eqnarray}}
\newcommand{\eeqn}{\end{eqnarray}\par\noindent}

\title{Compositionality as we see it, everywhere around us}  
\author{Bob Coecke\\ 
Cambridge Quantum Ltd.\\
Compositional Intelligence team, Oxford.\\
{\tt bob.coecke@cambridgequantum.com}} 

\begin{document}   
\maketitle  

\begin{abstract}
There are different meanings of the term ``compositionality" within science: what one researcher would call compositional, is not at all compositional for another researcher. The most established conception is usually  attributed to Frege, and is characterised by a bottom-up flow of meanings: the meaning of the whole can be derived from the meanings of the parts, and how these parts are structured together.  

Inspired by work on compositionality in quantum theory, and categorical quantum mechanics in particular, we propose the notions of Schr\"odinger,  Whitehead, and complete compositionality.  Accounting for recent important developments in quantum technology and artificial intelligence, these do not have the bottom-up meaning flow as part of their definitions.

Schr\"odinger compositionality  accommodates quantum theory, and also meaning-as-context. 
Complete compositionality further strengthens  Schr\"odinger compositionality in order to single out theories like ZX-calculus, that are complete with regard to the intended model.  All together, our new notions aim to capture the fact that compositionality is at its best when it is `real', `non-trivial', and even more when it also is `complete'.  

At this point we only put forward the intuitive and/or restricted formal definitions, and leave a fully comprehensive definition to future collaborative work.
\end{abstract}   

\section{Introduction}   

I was invited to write this paper for a special issue, aimed to celebrate Andrei Khrennikov's 60 birthday in recognition of his contributions to the quantum foundations community,  as well as Andrei being one of the frontrunners in using quantum structures in areas other than physics -- see Section \ref{sec:pic} for a picture of us.  In particular, we were asked to discuss why it makes sense to use certain quantum structures outside of physics, and in particular, why quantum physics provides structures that are useful outside of physics -- a question initially posed by Anton Zeilinger.   

In my work, already for some 20 years, the quantum structure that has taken the centre stage is the quantum compositional structure.  By this we refer to the particular manner in which quantum systems compose, which Schr\"odinger singled out as the key characteristic of quantum theory.  All-together, we provided a novel formalism for quantum theory entirely based on composition of systems, casting orthonormal bases \cite{CPV}, measurements \cite{CPav}, classicality \cite{CPaqPav}, and (strong) complementarity \cite{CD1, CD2, CDKZ},  entirely in compositional terms. The formalism can now be used to give a fully compositional account of all that  what one expects to get in a course on quantum foundations and quantum computing \cite{CKbook}.  It  is also now being used in quantum industry, e.g.~ZX-calculus based circuit optimisation techniques such as  \cite{clifford-simp, KissingerTcount, de2020fast, backens2020there} are now internal to Cambridge Quantum's t$|$ket$\rangle$ \cite{sivarajah2020t}.  Other uses in industry include \cite{Gidney2019}, and there are also efforts by major quantum computing companies to push our compositional quantum formalism as the number one educational tool for quantum computing. The key reason for that is that compositionality yields an entirely pictorial quantum formalism, usually referred to as `quantum picturalism' \cite{ContPhys}.

There of course have been earlier attempts to provide quantum theory with a novel formalism, see e.g.~\cite{CMW} for an overview.  While this quest was initiated by the great John von Neumann \cite{vN, BvN}, it is fair to say that none of those approaches have been successful -- with the exception of some applications outside of physics, e.g.~in NLP \cite{Widdows, widdows2003word}. One of the major issues of these attempts was their failure to account for quantum composition.  Meanwhile, many approaches in quantum foundations now also recognise the importance of the compositional backbone \cite{HardyJTF, Chiri1, CES, Chiri2, selby2017leaks, selby2018reconstructing, schmid2020unscrambling}. All of this seems to strongly confirm Schr\"odinger's hunch that compositionality underpins quantum theory.  

With respect to using quantum structures in other areas, in 2008 we used quantum compositional structure in order to solve a then open problem \cite{CCS, CSC}:  how can we combine meanings of words in order to produce the meaning of a sentence, with meanings represented as vectors, as is standardly done in modern natural language processing (NLP) \cite{Firth, harris1954distributional, PulmanDistributional}.  The crux was the fact that quantum compositional structure is, category-theoretically speaking, exactly the same as grammatical structure \cite{teleling}.  This correspondence has recently resulted in implementations of NLP on existing quantum computers, by ourselves at Cambridge Quantum \cite{QNLP-foundations, qspeak, Nature, QNLPPlus100, kartsaklis2021lambeq}. We provide more details on this in Sections \ref{sec:NLP} and  \ref{sec:top-down}.     

Now, speaking about different disciplines, the term ``compositionality" has different incarnations within different communities, and even within single communities there have been contradictory conceptions.  In the  recently emerging and flourishing  Applied Category Theory community, which considers very similar compositional structures as the ones discussed above, there does seem to be a  common understanding of what compositionality means.  In many cases it also manifests itself in terms of a pictorial formalism. Within this community, compositionality  has become a discipline transcending vehicle, that provides a bridge for exchanging concepts, methods and results.   Within one discipline it can provide a radically new perspective, putting the focus on how things interact, rather than what they are made up from.   However, as far as we know, nobody within this community has ever attempted to define  what the particular brands of  compositionality in use are.  This paper is our attempt to do so.

So in mathematical terms, compositionality is closely related to category theory.  However, there are some caveats.  Firstly, category theory was proposed as a meta-theory of mathematics, and the CT community at large pushed this view, while for us compositionality is all about real-world phenomena.  Secondly, many subfields of category theory favour so-called cartesian categories, while compositionally speaking these are the much less interesting categories. 

The multidisciplinary community mentioned above has its own journal and conference, although with different names, Compositionality and Applied Category Theory respectively.  Why? This dual terminology is a result of a compromise between communities with very different cultures. 
We favour ``Compositionality" as it clearly emphasises the key foundational concept that governs the work, while ``applied [something]" somewhat goes against the spirit of being of a foundational nature.  
On the other hand, in defence of the term ``Applied Category Theory", it refers to a clear community shift within the new generation of category-theoreticians, not just focusing on mathematics anymore.  This shift has many different dimensions to it,  one being the desire to solve important real-world problems. 

The oldest area doing practical ACT is computer science, thinking of data-types as systems, and programs as processes \cite{AbrLNP}. There is a huge literature, and every year  conferences produce a humongous amount of new papers -- some more history on compositionality in computer science is in \cite{DuskoTwitter}.  Other papers giving nice overviews of compositional structures in a variety of fields include \cite{BaezLNP, fong2018seven}, and a number of recent papers showing the wide range of topics currently addressed are \cite{ghica2007geometry, pavlovic2013monoidal, SobocinskiSignal, BaezFongElec, hinze2016equational, ConcSpacI, ghani2018compositional, smithe2020bayesian, shiebler2021category}. 

\section{Compositionality in linguistics and quantum physics} 


The area mostly associated with a formal definition of compositionality is formal linguistics, where its definition is usually stated as follows: 
\bit
\item The meaning of a whole (cf.~sentence) should only depend on the meanings of its parts (cf.~words) and how they are fitted together (cf.~grammar).  
\eit
A prominent argument for this notion, put forward by Frege \cite{frege1914letter}, is that without compositionality it would be impossible for us to understand a sentence that we have never heard before.  Broadly speaking, this is of course indeed the case.  However, interestingly Frege himself mentions in `Grundlagen der Arithmetik' \cite{frege1884grundlagen, janssen2001frege}: 
\bit
\item ``Never ask for the meaning of a word in isolation, but only in the context of a sentence."
\eit
This is  called `the context principle' \cite{PhilosophicalInvestigations}, and indeed: 
\bit
\item In language, meanings are to a great extent induced by the context. For example, ambiguous words such as {\tt Alice}  could point at either of these two incarnations of {\tt Alice}:
\[
\epsfig{figure=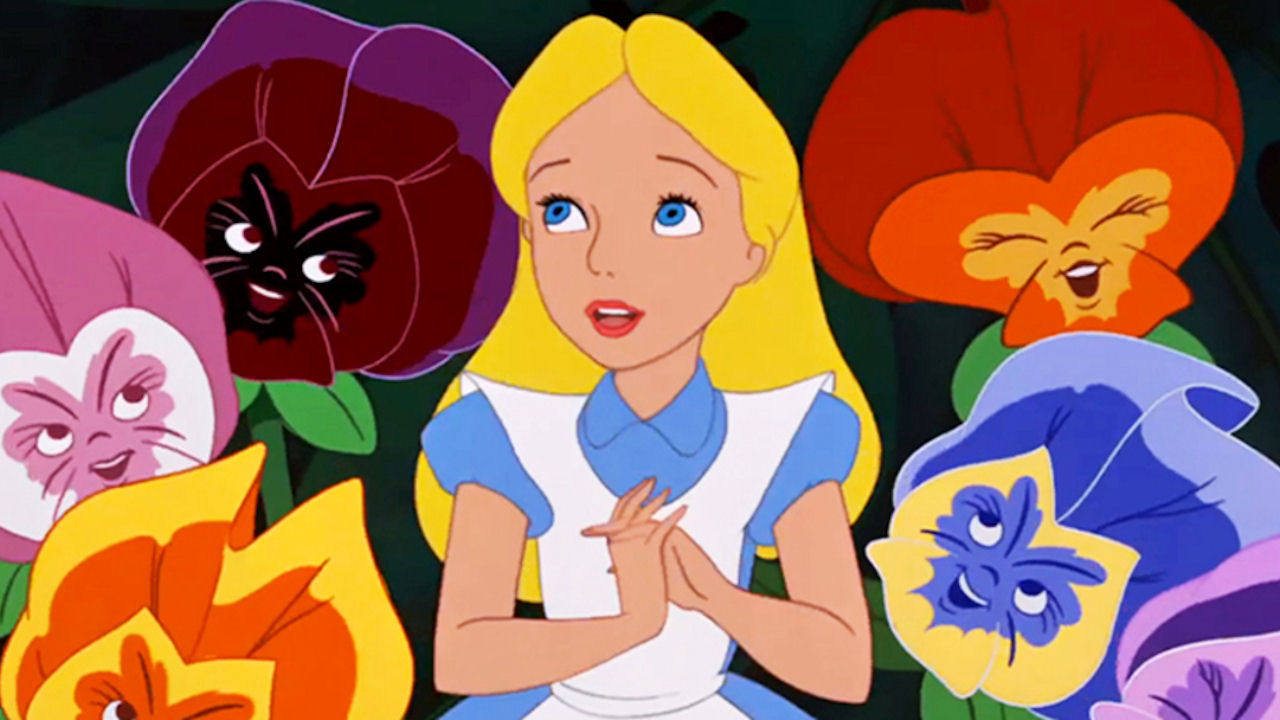,height=80pt}\qquad\ \  \qquad\qquad\epsfig{figure=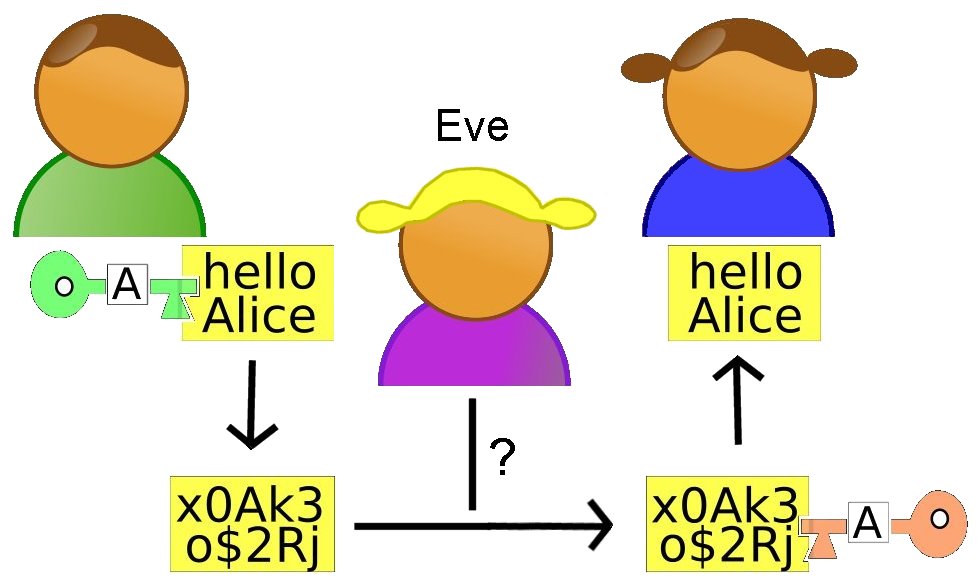,height=80pt} 
\]
We may disambiguate {\tt Alice}  by other words in the sentence or  larger text, or by non-linguistic features, like the conversation taking place between two children impersonating fairytale characters, or between security protocol experts. Even within the restricted context of gothic rock fans {\tt Alice} could still  have two meanings, two cases being:    
\[
\epsfig{figure=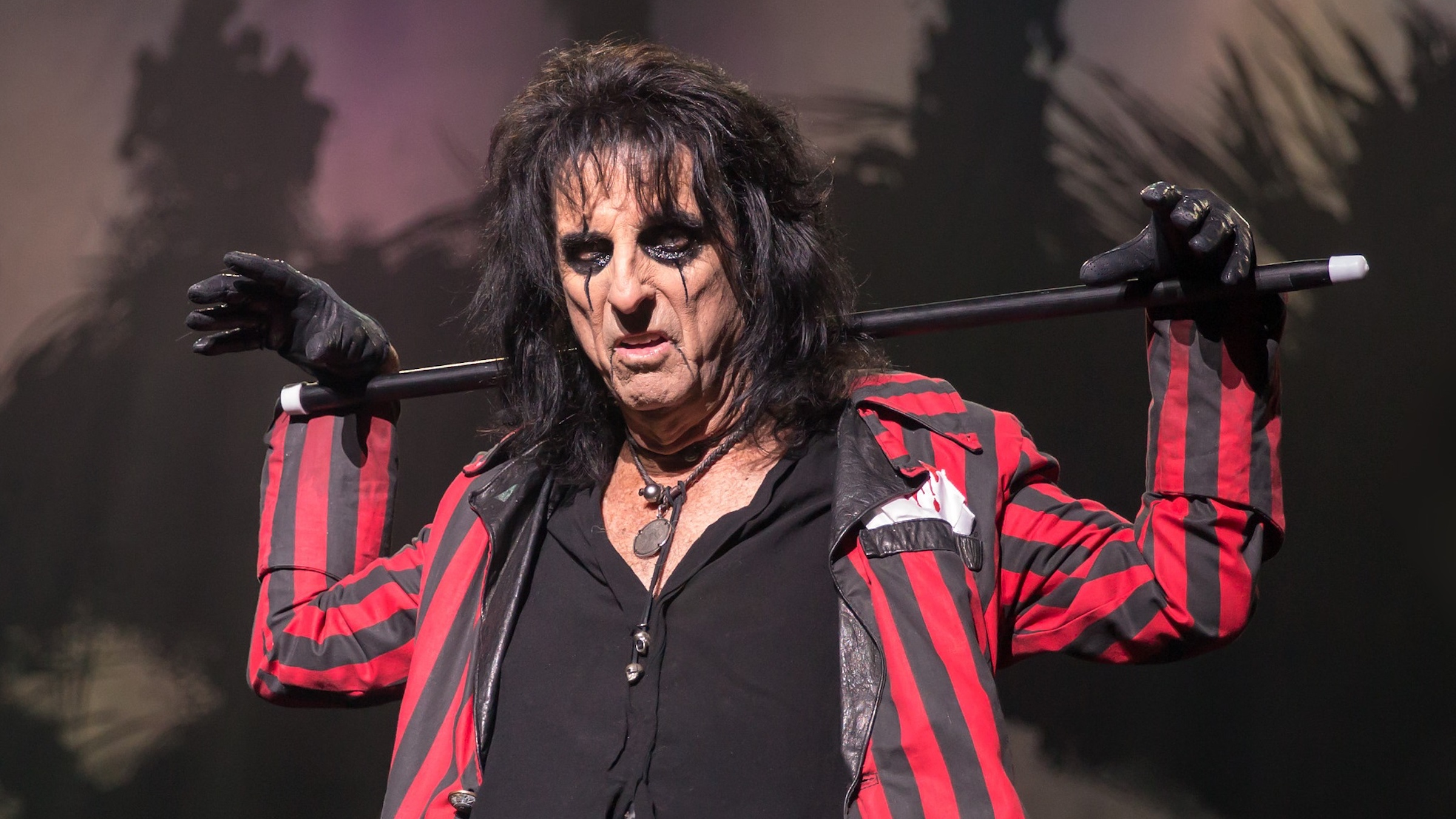,height=80pt}\qquad\ \quad  \qquad\qquad\epsfig{figure=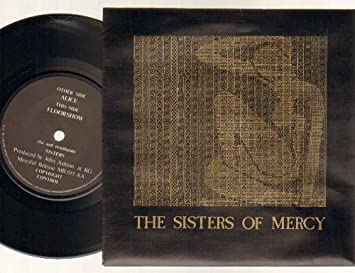,height=80pt}  \ \ \ \ 
\]
\eit
There are some more refined formulations of linguistic compositionality that aim to account for context-dependent word ambiguities, however, the following is harder to get around:
\bit
\item In language,  separation between word meanings and how they grammatically fit together doesn't always make sense, as the grammar can  depend on the  overall meaning, which may  need   additional context to be disambiguated. For example, in {\tt black metal fan} the grammar depends on the meanings:\footnote{We took this particular example  from \cite{CoeckeMeich}, but there are of course many more.}
\[
\ \ \ \quad\epsfig{figure=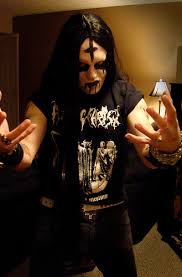,height=140pt}\qquad\ \  \qquad\qquad\epsfig{figure=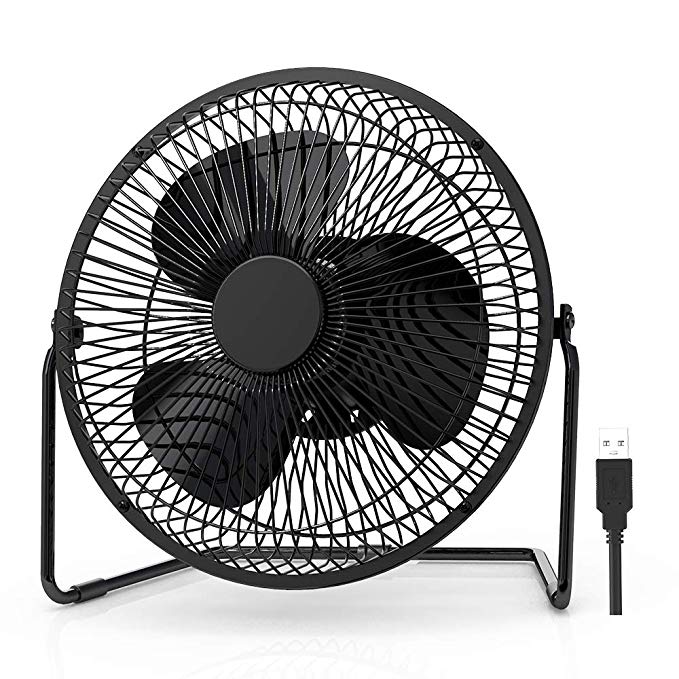,height=140pt} 
\]
Using the language of \cite{CoeckeText}, the respective grammatical parsings proofs are: 
\[
\input{Fan3.tikz}\qquad\qquad\qquad\input{Fan2.tikz} 
\]
It is not as important to understand them, than just to observe that they are different, so there is no fixed grammar until the overall meaning has been unambiguously established.
\eit
In the linguistic context the fundamental problem with the above stated notion of compositionality is that it puts a special focus on bottom-up meaning-flows, the meaning of the whole being induced by the meanings of the parts.  However, meaning-flows within a composite can equally well travel top-down along the same structure, that is, the meaning of the whole can contribute to meanings of parts. 
So what we are aiming for is some network structure that allows for such two-way meaning-flows.  Maybe even better put, a network structure that allows for meaning-relationships along the different scales, cf.~word vs.~sentence vs.~text.
 
In quantum theory the situation gets even more extreme:
\bit
\item  For entangled states the idea that subsystems have meanings breaks down, and so does Frege compositionality.  For example, the pure Bell-state's sub-systems are in the state of total noise when conceived individually. In addition to these individual descriptions,   in order to know that the overall state is the Bell-state, one needs the entire description of the Bell-state itself.  
\eit   
So  there is no meaningful decomposition of the whole in terms of parts plus structure. Yet both language and quantum theory are perfectly compositional as far as we are concerned. How?
 
\section{Process-theoretic compositionality} 
 
Let's take quantum theory as our starting point and see how we can associate a notion of compositionality to it that  generalises to language. As already pointed out in the introduction, Schr\"odinger singled out the  nature of composition of systems  in quantum theory not as just a, but as \underline{the} defining feature of the theory \cite{SchrodingerComp}:
\begin{quote}
``Entanglement is not one but rather the characteristic trait of quantum mechanics."
\end{quote}
Formally, (pure-state) entanglement means:
\beq\label{eq:sep}
|\psi_{AB}\rangle \not= |\psi_{A}\rangle\otimes |\psi_{B}\rangle
\eeq 
that is, the state of two systems $A$ and $B$ doesn't factor as a state  $|\psi_{A}\rangle$ of $A$ and a state  $|\psi_{B}\rangle$ of $B$.  
We can write (\ref{eq:sep}) diagrammatically as follows \cite{Kindergarten, CKbook}:
\beq\label{eq:sep2}
\input{separate.tikz}
\eeq
This notation of states is a special case of representing processes as boxes:
\[
\input{box2.tikz}
\]
which can be composed to form \em diagrams\em: 
\[
\input{boxexample2.tikz}
\]
An example of a maximally  entangled state/effect are the \em Bell-state/effect\em, which we denote as:
\[
\input{cap.tikz}\ = \ket{00} + \ket{11}\qquad\qquad\qquad\input{cup.tikz} \ = \bra{00} + \bra{11}
\]
With this state and effect our diagram becomes:
\[
\input{boxexample3.tikz}
\]

The claim we made earlier that the  Bell-state's sub-systems are in the state of  noise when conceived individually, diagrammatically takes the form:
\[
\input{capdiscard.tikz} 
\]
We refer to a situation where the composite of two systems cannot be meaningfully described  in terms of its parts as \em non-trivial  composition\em.  Others may call it wholistic.  The Bell-states, as well as any other maximally entangled states,  can be shown to constitute an extremal example of non-trivial composition \cite{coecke2018uniqueness, coecke2019uniqueness}.  Formally, the possible states of two systems are then isomorphically related to the processes that each system can undergo as follows:
\[
\left\{ \input{box5.tikz} \right\}\ \simeq\ \left\{ \input{box4.tikz} \right\} 
\]
where the isomorphism is established by the following  equation:
\[
\input{yank.tikz} 
\]
which is characteristic for extremal  non-trivial composition. 

What makes these kinds of states (and effects) interesting is not what they are, but what one can do with them.  That is, what happens when one processes them.  This for example leads to features like quantum teleportation \cite{Tele, swap, TeleExp},  quantum non-locality \cite{Bell, GHZ} and, measurement-based quantum computing \cite{Gottesman, MBQC2}, among   many more, each of these having no non-trivial classical counterparts.  

This means that a general notion of  compositionality should apply to general theories of processes, not just theories of states. With the advent of our compositional quantum formalism, a.k.a.~categorical quantum mechanics \cite{AC1, CKbook}, a.k.a.~quantum picturalism \cite{ContPhys}, we know that such a presentation leads to a much more intuitive purely diagrammatic formalism, important computational advantages, e.g.~in quantum computing \cite{DP2, de2017zx, duncan2019graph, kissinger2019pyzx}, and bridges disciplines giving rise to new technologies such as quantum natural language processing \cite{WillC, QNLP-foundations, QNLPPlus100}. The paper \cite{coecke2021kindergarden} lists  recent compositionality-driven advances for the specific  case of the ZX-calculus.    

The idea of focusing on processes instead of states is of course not new, and even goes back to the pre-Socatic Heraclitus' ``panta rhei", a.k.a.~everything flows, and in  more recent times was  advocated by Whitehead  in ``Process and Reality" \cite{Whitehead}. His stance departs from our typical Western metaphysics in which we start with a static kinematic layer with dynamics as a secondary derivative \cite{seibt2012process}.  Whitehead's work however did not explicitly consider a non-trivial \em parallel \em composition as in  (\ref{eq:sep2}), and therefore we will introduce two notions of compositionality, distinguishing between Whitehead and Schr\"odinger compositionality.

Diagrammatically, we can think of Whitehead compositionality as diagrams like this: 
\beq\label{eq:box3}
\input{box3.tikz} 
\eeq
Condition (\ref{eq:sep2}) for two systems $A$ and $B$ `separated in space', has an analogue for processes with input $A$ and output $B$ `separated in time':
\beq\label{eq:sep3}
\input{separate2.tikz} 
\eeq
This is satisfied for any reasonable theory for processes, since otherwise we would obtain a theory only containing constant processes that produce output $|\psi_{B}\rangle$ for any possible input: 
\[
f(|\phi\rangle) = |\psi_{B}\rangle\langle \psi_{A}|\phi\rangle \simeq |\psi_{B}\rangle
\]
Similarly, processes also enjoy non-triviality by analogy, given that what they typically do is  relate varying inputs to corresponding outputs.  Hence deriving them from a fixed input-output pair doesn't make much sense.

For illustrative purposes we will mainly make use of the conception of process theory as in \cite{CKpaperI, CKbook}. It has meanwhile been established that the process theory format provides a useful complementary presentation of quantum theory as compared the the usual Hilbert space formulation. In particular, in the foundations of quantum theory a process-theoretic backbone is now a dominant feature for many developments -- see e.g.~\cite{HardyJTF, Chiri1, CES, Chiri2, selby2017leaks, selby2018reconstructing, schmid2020unscrambling}.   

\begin{defn}{\cite{CKpaperI, CKbook}}\label{def:ProcTheor}
A \em process theory \em is defined to consist of:
\bit
\item a collection of systems (or types, or interfaces) $S$  represented by wires,  
\item a collection of processes $P$ represented by  boxes, and where each process in $P$ has a number of input wires and a number of output wires taken from $S$, and  
\item a means of composing processes, represented by wiring boxes together, and the result of doing so should again be a process in $P$.
\eit
\end{defn}

While in quantum theory one can interpret these processes as taking place in space and time, this doesn't always have to be the case, for example, wirings can be physical wires:
\[
   \input{electricity-example2.tikz} \quad := \quad
  \raisebox{-2.7cm}{\includegraphics[width=2.4cm]{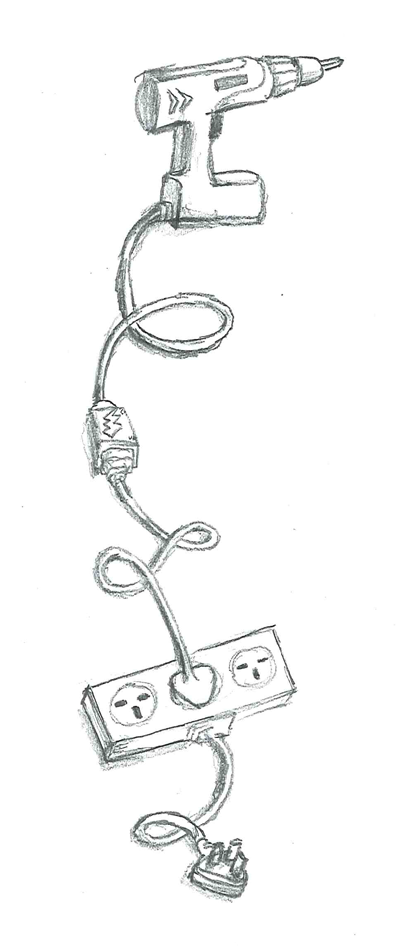}} 
\]
or, more abstract connections like in a flow-chart \cite{FlowChart}. 
  
However, the definition of process theories stated above is by no means the most general notion conceivable of a compositional theory.  For example, besides processes a theory could also specify transformations of processes (a.k.a.~higher-order processes---see e.g.~\cite{Chiribella2013, Oreshkov2012, kissinger2017categorical}), and it could also have transformations of interfaces etc.  In the case of the latter, one could think of these diagrammatically as surfaces connecting wires, and processes transforming processes are typically denoted as `combs' \cite{comb}:  
\[
\input{comb.tikz}
\]
We could then also have hyper-combs transforming combs, and so on.  
The bottom line is that our definitions below also apply to these more general notions of compositional theories that stretch well beyond the notion of process theory stated above.

Some of the above can also be framed in the language of category theory.  In the case of process theories \`a la  Whitehead we would be talking about ordinary categories, possibly cartesian, while in the case of process theories \`a la Schr\"odinger we would be talking about proper  (strict) symmetric monoidal categories  \cite{MacLane, CatsII}. In the extremal case of Bell-states we are talking about compact closed categories \cite{Kelly, KellyLaplaza, AC1, CatsII}. Going beyond process theories  by considering transformations of processes and so on takes us into higher category theory \cite{leinsterbook}. 

\section{Kinds of compositionality}

\subsection{Schr\"odinger compositionality}   

\begin{defn}\label{defn:schro}
A \em Schr\"odinger compositional theory \em is a compositional theory in the sense of the previous section, subject to two conditions: 
\bit
\item  Composition is non-trivial, that is, the description of a whole cannot always be decomposed in any meaningful way---cf.~given the description of parts, one may still require a full description of the whole in order to derive the latter. 
\item The ingredients of a compositional theory---e.g.~the ingredients in Defn.~\ref{def:ProcTheor}---all have clear meaningful ontological counterparts in reality.
\eit
\end{defn}

The first condition of non-triviality is what makes the passage to a compositional description worthwhile.  If it is trivial (cf.~cartesian categories), then there really is not much point to elevate composition to a 1st class citizen.  In fact, it is our belief that pretty much everything in the real world admits non-trivial composition when regarded in process-theoretic terms.  

The second condition of the compositional structure `being real' is essential in that it provides a model for the theory, hence establishing its consistency, and much more importantly, makes the theory relevant.  It also provides us with our daily intuition as guidance for studying it, in the same way that Euclid relied on physical space to establish geometry.  Finally, for the specific case of AI, compositionality then becomes a vehicle for incorporating embodiment in a compositional fashion, as for example in \cite{ConcSpacI, TalkSpace}.

So a  slogan could be:
\begin{center}
\fbox{\em Compositionality is at its best when it is non-trivial and real.} 
\end {center}

The main purpose of the following definition is  just to contrast it with Defn.~\ref{defn:schro}: 

\begin{defn}
A \em Whitehead compositional theory \em is the same as a Schr\"odinger compositional theory except for the fact that composition is only one-dimensional---cf.~(\ref{eq:box3}).
\end{defn}

\subsection{Complete compositionality}

One issue with Schr\"odinger compositionality as we defined it is that we assumed processes to be `black-boxes', that is, they don't necessarily interact in any obvious way.  Ultimately, this boils down to what the equational structure of the process theory is: some may have a lot of equations between diagrams, other few, or none at all if they are free theories.  However, there are some very special kinds of process theories around with a very rich compositional structure.   They typically have a number of generating processes, which have a very easily understood interaction structure, and also generate a very rich class of processes by plugging them together. This results in a very powerful reasoning system.   

 One example of such a rich compositional theory is the ZX-calculus \cite{CD1, CD2}, where, for example, a CNOT-gate arises from two generating processes \cite{coecke2021kindergarden}:
\beq\label{CNOT}
  \input{cnot.tikz} \ \ =\ \ \input{cnotfromspiders.tikz}  
\eeq
These generating processes respectively correspond to Z- and X-copying, or more generally, the respective  `spiders' they generate:
\[
\begin{array}{l}
\input{greenspider.tikz}\ \ :=\ \ \ket{0\ldots 0}\bra{0\ldots 0} + e^{i \alpha} \ket{1\ldots 1}\bra{1\ldots 1}
\vspace{4mm}\\
\input{redspider.tikz}\ \ :=\ \  \ket{+\ldots +}\bra{+\ldots +} + e^{i \alpha} \ket{-\ldots -}\bra{-\ldots -} 
\end{array}   
\]
These spiders have an interaction structure that can be written down in only 4 equations, but at the same time generate all matrices over arrays of qubits \cite{CKbook}.
Another example of a complete compositional theory is ZW-calculus \cite{CK, Amar}, for which the generating processes are the GHZ-state and the W-state.   

How can we qualify the richness of such a theory?  We can do this in terms of how much of all the intended equations between diagrams can be derived from the interaction rules for the generating processes.  In the case of the ZX-calculus and the ZW-calculus, in fact, all equations that can be derived using linear algebra can be derived from the interaction rules \cite{hadzihasanovic2018two}. We say that both of these calculi are \em complete\em. This completeness caused quite some noise at the time, but it didn't come entirely out of the blue, some earlier work including \cite{Backens, Backens2, Amar, jeandel2018complete}, and there also are some further elaborations \cite{DBLP:conf/rc/CoeckeW18, vilmart2019near}.    
 
\begin{defn}
A \em complete compositional theory \em is a Schr\"odinger compositional theory with a complete compositional structure in the above sense.
\end{defn}

Our slogan now becomes:
\begin{center}
\fbox{\em Compositionality is at its best when it is non-trivial, real  and complete.} 
\end {center}

There is one subtlety however...  

\subsection{Relativity of compositionality}

In the case of ZX-calculus by `real' we mean that each spider corresponds to a physical process that may occur.  We say `may', as most spiders cannot be realised in a deterministic fashion: they can only occur as a branch in a non-deterministic process, e.g.~an appropriately crafted measurement.  In other words, they are not `causal' \cite{Chiri1, coecke2016terminality, kissinger2017equivalence}.  

Now, rather than considering all quantum processes, one could restrict to unitary ones.  In that context, the ZX-spiders are not `real' anymore and only occur as `1/2 of a CNOT-gate', as in (\ref{CNOT}).  So now  real processes are further decomposed in smaller bits.  These smaller bits greatly facilitate computation as compared to the unitary processes.  

\begin{defn}
\em Lego compositionality \em is a genealisation of Schr\"odinger compositionality in that  processes may be broken down in smaller pieces that have no counterpart in reality.  
\end{defn}

\section{Examples and non-examples}  

\subsection{Natural language processing}\label{sec:NLP}

Our approach to natural language from  \cite{CCS, CSC, GrefSadr, KartSadr} called DisCoCat, which combines meaning and linguistic structure, is Schr\"odinger compositional, and has equally non-trivial composition  as quantum theory.  Boxes in diagrams such as:
\beq\label{eq:flowers1rp}
\input{flowers1rp.tikz} 
\eeq
correspond to word-meanings, the wirings correspond to grammatical structure, and the resulting state is the meaning of the sentence. Hence the compositional structure is clearly `real'.  

The word-meanings themselves have composite types, and a transitive verb box such as: 
\[
\input{flowers1rp2.tikz} 
\]
evidently satisfies (\ref{eq:sep2}) since otherwise the subject and the object would play no role in the overall meaning of the sentence.  An adjective like:\footnote{This example was 1st pointed out to us by Fabrizio Genovese.}
\[
\input{fake.tikz} 
\]
cannot really be assigned any meaning until we actually know what its context is, so  there clearly is no meaningful decomposition of it over parts and structure.
 
The wiring of {\tt that} in terms of spiders taken from \cite{FrobMeanI} can be considered as complete compositional, given that a single relative pronoun is broken down in smaller bits -- formally, when interpreted in vector space, a GHZ-state and a $|+\rangle$-state.  We recently proposed many more such decompositions in \cite{GramEqs}, and there is a case to be made that they are in fact `real', as they result in new truths about grammar not previously recorded.  

These language diagrams can be equivalently interpreted as quantum physics diagrams \cite{teleling, QNLP-foundations}, a connection that can only be established when conceiving both theories in the compositional manner we casted them in.  This is an important   feature about  of shared compositional structure: it enables one to map seemingly entirely different  areas onto each other, and enables an exchange of concepts and tools.  In fact, in this particular case it gets even better, as this connection effectively enabled us to process natural language on a quantum computer \cite{WillC, QPL-QNLP, Nature, QNLP-foundations, QNLPPlus100, kartsaklis2021lambeq}.  
For doing so, a language diagram, now interpreted as a quantum process, still needs to be compiled so that it fits on quantum hardware.  Again this compilation is enabled by compositional structure, in this case the complete compositionality of the ZX-calculus, by `filling' the meaning-boxes with ZX-diagrams, which can then be put in circuit-form \cite{QNLP-foundations}:
\beq\label{eq:abs} 
\input{abs1.tikz}\ \ \ \mapsto\ \ \ \input{abs2.tikz}  
\eeq

\subsection{Neural networks and tensor networks} 

In neural networks, individual neurons typically have no clear counterpart in reality, and due to how the non-linearities are inserted, have no comprehensible interaction rule.  Hence standard neural networks do not satisfy the definitions of compositionality put forward here. 

In the case of tensor networks in physics, nodes can have counterparts in condensed matter, but in many cases again there is no rich structure of interaction rules.  Hence, tensor networks can be Schr\"odinger compositional, but are typically nowhere near complete compositional.

\subsection{Can maths be Schr\"odinger compositional?}

Well, it kind of depends whether one conceives maths as either real or not, so it is to some extent subjective.  Something that one can do is to recast mathematical structures in Schr\"odinger compositional terms. We for example did this with convex structures in \cite{ConcSpacI}, which following G\"ardenfors \cite{gardenfors, gardenfors2} can be interpreted as related to our sensual perception, and our theory then expresses how different senses compose, generating more sophisticated concepts from simple concepts.  Similarly, in \cite{TalkSpace} we take 3(+1)D cartesian space, and recast it as compositional spatial relations that represent how entities within physical space are related.  
 
\section{Top-down flow and ambiguity}\label{sec:top-down}  

Formal linguistics  is the area most associated with compositionality.  But, categorial grammars such as in \cite{Ajdukiewicz, Bar-Hillel, Lambek0, montague1970universal, montague1973proper, lambekbook} have a  surprising limitation   in that compositionality doesn't stretch beyond the sentence level, while, of course, when we speak/write we compose sentences in order to form larger chunks of text. DisCoCat, as described above,  having been built on top of  categorial grammars and is therefore also limited to relating meanings of words to those of sentences, and not to larger text.

Our attempt to fix this shortcoming was started in \cite{CoeckeText} under the name DisCoCirc  and is being pushed further in \cite{GramCircs}.   The main slogan of our work is:
\begin{center}
\fbox{\em A sentence is a process that alters the meanings of its inhabitants.} 
\end {center}
The idea is to move from diagrams like (\ref{eq:flowers1rp}) to circuits:  
\[
\input{flowersrpQUINT3rev.tikz} 
\]
which can then be composed further: 
\[
\input{flowersrpQUINT3revcopy.tikz} 
\]
A major feature about these \emph{language circuits} is that they can be interpreted bottom-up, as well as top-down, or a hybrid that expresses relationships between different levels.  For example, we could have the meanings of a circuit as a given, and from these deduce the meanings of the words in it.  This is in fact what we did in our QNLP-work \cite{Nature, QNLP-foundations, QNLPPlus100, kartsaklis2021lambeq}, obtaining the word-meanings from a corpus of sentences we knew to be truthful.  The reason for doing so in that particular case is that currently there is no better way available to load classical data on a quantum computer. So this is a case where a top-down interpretation is imposed by the technology. 


As mentioned in the beginning of this paper, one of the key features of language is ambiguity, which manifests in many different ways. Evidently, when reasoning top-down, there initially will be ambiguity about the parts. In modern natural language processing meanings are represented by vectors.   But wait a minute, we know how to deal with ambiguity in quantum theory: mixed states.  We can now use our compositional correspondence to transfer von Neumann's quantum account on ambiguity in terms of density matrices to language \cite{calco2015}:
\beqa
\rho_{\mbox{\scriptsize\tt queen}} &=& |\mbox{\tt queen-royal}\rangle\langle\mbox{\tt queen-royal}|   \\
&+&|\mbox{\tt queen-band}\rangle\langle\mbox{\tt queen-band}|   \\
&+&|\mbox{\tt queen-bee}\rangle\langle\mbox{\tt queen-bee}|   \\
&+&|\mbox{\tt queen-chess}\rangle\langle\mbox{\tt queen-chess}|   \\
&+&|\mbox{\tt queen-card}\rangle\langle\mbox{\tt queen-card}| \\
&+& \ldots
\eeqa
This can then be extended to all kinds of ambiguities, including different grammars  with different meanings like in the case of {\tt black metal fan}:  
\beq\label{mixture}
\left|\mbox{\raisebox{-14pt}{\,\epsfig{figure=FanPerson.jpg,height=32pt}\,}}\right\rangle\left\langle\mbox{\raisebox{-14pt}{\,\epsfig{figure=FanPerson.jpg,height=32pt}\,}}\right| + \left|\mbox{\raisebox{-14pt}{\,\epsfig{figure=FanDevice.jpg,height=32pt}\,}}\right\rangle\left\langle\mbox{\raisebox{-14pt}{\,\epsfig{figure=FanDevice.jpg,height=32pt}\,}}\right|
\eeq  
The top-down meaning-flow will then disambiguate such mixtures. 

You may think that the sum in (\ref{mixture}) takes us out of the  fully compositional paradigm, but  in 2005 Selinger put forward a paradigm that generates mixedness just using diagrams \cite{SelingerCPM}, used throughout our book \cite{CKbook}. It's where the `ground' notation also entered the quantum formalism \cite{CPer}:
\[
\input{ground.tikz} 
\]
It has even recently obtained the status of being completely compositional in  \cite{carette2019completeness}. 


\section{Final considerations}

\paragraph{Why quantum?} So why did the notion of Schr\"odinger compositionality emerge from quantum theory? No particular reason, really.  Historically, of course,  the Western metaphysical bottom-up perspective favoured reductionist interpretations. What quantum added to the picture was that it puts a case forward  that simply couldn't be ignored.  
We also strongly believe that this is the kind of math that is very needed in humanities and social sciences, so that we finally can relieve ourselves of the mindlessness of mostly doing statistics. Ultimately, we think of these compositional structures as everywhere around us, a case that we have made at a number of previous occasions \cite{ContPhys, DBLP:books/daglib/p/Coecke17, qspeak}.     
 
\paragraph{What's next?} Admittedly, we kept things informal, adopting a restricted notion of processes,  and stated the conditions on Schr\"odinger and complete compositionality only informally. We are currently working towards a fully-formal characterisation of several notions of compositionality.  

\section{Pic with Andrei}\label{sec:pic}
 
One of my first better-known papers in the area of compositionality was the transcription of a talk at Andrei Khrennikov's annual quantum foundations conference, in 2005, called Kindergarten Quantum Mechanics \cite{Kindergarten}.  The next year, in 2006,  Andrei got us in the local Swedish newspaper after I gave a talk ``In the beginning God created $\otimes$'', translated as ``I begynnelsen skapade Gud sp\"annmuskeln som en bild": 
\[
\includegraphics[width=14cm]{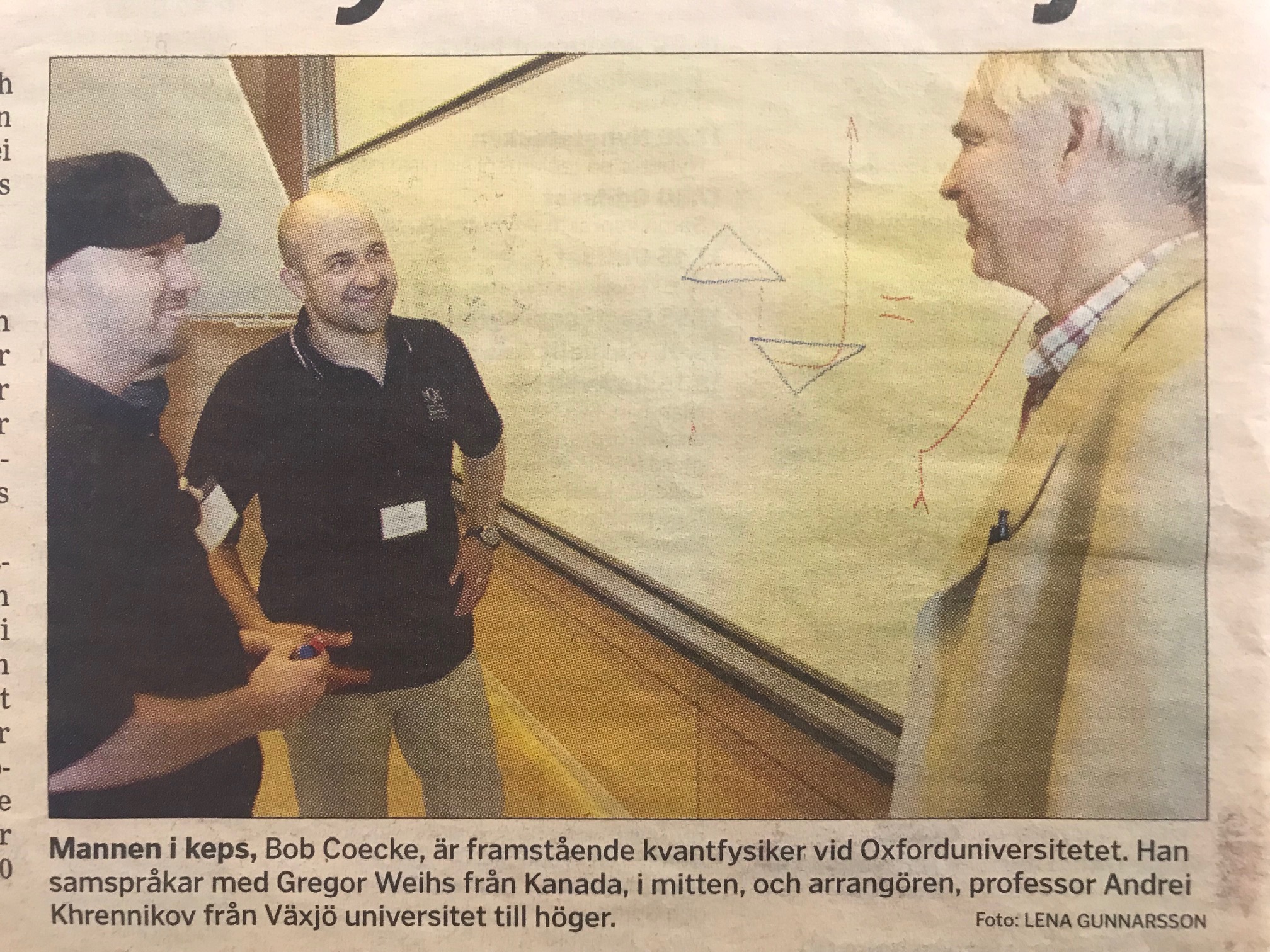} 
\]
We thank Vincent Wang for proofreading this manuscript.  At the time that this picture was taken, Vincent was still playing with teddy bears. Stephen Clark, Amar Hadzihasanovic,  
Arkady Plotnitsky,  Ilyas Khan and Ian Durham also provided feed-back.  Peter Stannack pointed out that this paper can be conceived as providing a response to a puzzle put forward by Ryan Nefdt in \cite{nefdt2020puzzle}. 
 
%
%
%
%
%
%
%

\bibliographystyle{plain} 
\bibliography{mainNOW}

\end{document}

%% file: macros/defs.tex
\newkeycommand{\boxmap}[style=small box][1]{\,\tikz{\node[style=\commandkey{style}] (x) {$#1$};\draw(0,-1.25)--(x)--(0,1.25);}\,}
\newkeycommand{\dboxmap}[style=dbox][1]{\,\tikz{\node[style=\commandkey{style}] (x) {$#1$};\draw[boldedge] (0,-1.25)--(x)--(0,1.25);}\,}

\newkeycommand{\wirelabel}[style=white label]{\,\tikz{\node[style=\commandkey{style}]{};}\,\xspace}

\newkeycommand{\pointmap}[style=point][1]{\,\tikz{\node[style=\commandkey{style}] (x) at (0,-0.3) {$#1$};\node [style=none] (3) at (0, 0.9) {};\node [style=none] (3) at (0, -0.9) {};\draw (x)--(0,0.7);}\,}

\newkeycommand{\point}[style=point,estyle=][1]{\,\tikz{\node[style=\commandkey{style}] (x) at (0,-0.05) {$#1$};\draw [\commandkey{estyle}] (x)--(0,1.05);}\,}

\newkeycommand{\copointmap}[style=copoint][1]{\,\tikz{\node[style=\commandkey{style}] (x) at (0,0.3) {$#1$};\draw (x)--(0,-0.7);}\,}

\newkeycommand{\dpoint}[style=point][1]{\,\tikz{\node[style=\commandkey{style},doubled] (x) at (0,-0.3) {$#1$};\draw[boldedge] (x)--(0,0.7);}\,}

\newkeycommand{\kpoint}[style=kpoint,estyle=][1]{\,\tikz{\node[style=\commandkey{style}] (x) at (0,-0.05) {$#1$};\draw [\commandkey{estyle}] (x)--(0,1.05);}\,}

\newkeycommand{\kpointgray}[style=gray kpoint,estyle=][1]{\,\tikz{\node[style=\commandkey{style}] (x) at (0,-0.05) {$#1$};\draw [\commandkey{estyle}] (x)--(0,1.05);}\,}

\newkeycommand{\typedkpoint}[style=kpoint,edgestyle=][2]{%
\,\begin{tikzpicture}
  \begin{pgfonlayer}{nodelayer}
    \node [style=none] (0) at (0, 1) {};
    \node [style=\commandkey{style}] (1) at (0, -0.75) {$#2$};
    \node [style=right label] (2) at (0.25, 0.5) {$#1$};
  \end{pgfonlayer}
  \begin{pgfonlayer}{edgelayer}
    \draw [\commandkey{edgestyle}] (1) to (0.center);
  \end{pgfonlayer}
\end{tikzpicture}\,}

\newkeycommand{\typedkpointdag}[style=kpoint adjoint,edgestyle=][2]{%
\,\begin{tikzpicture}
  \begin{pgfonlayer}{nodelayer}
    \node [style=none] (0) at (0, -1) {};
    \node [style=\commandkey{style}] (1) at (0, 0.75) {$#2$};
    \node [style=right label] (2) at (0.25, -0.5) {$#1$};
  \end{pgfonlayer}
  \begin{pgfonlayer}{edgelayer}
    \draw [\commandkey{edgestyle}] (0.center) to (1);
  \end{pgfonlayer}
\end{tikzpicture}\,}

\newkeycommand{\bistate}[style=kpoint][1]{\,\begin{tikzpicture}
  \begin{pgfonlayer}{nodelayer}
    \node [style=\commandkey{style}, minimum width=1 cm, inner sep=2pt] (0) at (0, -0.25) {$#1$};
    \node [style=none] (1) at (-0.75, 0) {};
    \node [style=none] (2) at (0.75, 0) {};
    \node [style=none] (3) at (-0.75, 0.88) {};
    \node [style=none] (4) at (0.75, 0.88) {};
  \end{pgfonlayer}
  \begin{pgfonlayer}{edgelayer}
    \draw (1.center) to (3.center);
    \draw (2.center) to (4.center);
  \end{pgfonlayer}
\end{tikzpicture}\,}

\newkeycommand{\bistatebig}[style=kpoint][1]{\,\begin{tikzpicture}
  \begin{pgfonlayer}{nodelayer}
    \node [style=\commandkey{style}, minimum width=, minimum width=1.26 cm, inner sep=2pt] (0) at (0, -0.25) {$#1$};
    \node [style=none] (1) at (-1, 0) {};
    \node [style=none] (2) at (1, 0) {};
    \node [style=none] (3) at (-1, 0.88) {};
    \node [style=none] (4) at (1, 0.88) {};
  \end{pgfonlayer}
  \begin{pgfonlayer}{edgelayer}
    \draw (1.center) to (3.center);
    \draw (2.center) to (4.center);
  \end{pgfonlayer}
\end{tikzpicture}\,}

\newkeycommand{\dbistate}[style=dkpoint][1]{\,\begin{tikzpicture}
  \begin{pgfonlayer}{nodelayer}
    \node [style=\commandkey{style}, minimum width=1 cm, inner sep=2pt] (0) at (0, -0.25) {$#1$};
    \node [style=none] (1) at (-0.75, 0) {};
    \node [style=none] (2) at (0.75, 0) {};
    \node [style=none] (3) at (-0.75, 1.0) {};
    \node [style=none] (4) at (0.75, 1.0) {};
  \end{pgfonlayer}
  \begin{pgfonlayer}{edgelayer}
    \draw [boldedge] (1.center) to (3.center);
    \draw [boldedge] (2.center) to (4.center);
  \end{pgfonlayer}
\end{tikzpicture}\,}

\newkeycommand{\bistatebraket}[style1=kpoint,style2=kpointadj][2]{\,\begin{tikzpicture}
  \begin{pgfonlayer}{nodelayer}
    \node [style=\commandkey{style1}, minimum width=1 cm, inner sep=2pt] (0) at (0, -0.75) {$#1$};
    \node [style=none] (1) at (-0.75, -0.5) {};
    \node [style=none] (2) at (0.75, -0.5) {};
    \node [style=none] (3) at (-0.75, 0.5) {};
    \node [style=none] (4) at (0.75, 0.5) {};
    \node [style=\commandkey{style2}, minimum width=1 cm, inner sep=2pt] (5) at (0, 0.75) {$#2$};
  \end{pgfonlayer}
  \begin{pgfonlayer}{edgelayer}
    \draw (1.center) to (3.center);
    \draw (2.center) to (4.center);
  \end{pgfonlayer}
\end{tikzpicture}\,}

\newkeycommand{\weight}[style=dkpoint][1]{\,\begin{tikzpicture}
  \begin{pgfonlayer}{nodelayer}
    \node [style=\commandkey{style}] (0) at (0, -0.5) {$#1$};
    \node [style=upground] (1) at (0, 0.75) {};
  \end{pgfonlayer}
  \begin{pgfonlayer}{edgelayer}
    \draw [style=boldedge] (0) to (1);
  \end{pgfonlayer}
\end{tikzpicture}\,}

\newkeycommand{\dkpoint}[style=dkpoint][1]{\,\tikz{\node[style=\commandkey{style}] (x) at (0,-0.1) {$#1$};\draw[boldedge] (x)--(0,1);}\,}
\newkeycommand{\dkpointadj}[style=dkpointadj][1]{\,\tikz{\node[style=\commandkey{style}] (x) at (0,0.1) {$#1$};\draw[boldedge] (x)--(0,-1);}\,}
\newkeycommand{\dkpointtrans}[style=dkpointtrans][1]{\,\tikz{\node[style=\commandkey{style}] (x) at (0,0.1) {$#1$};\draw[boldedge] (x)--(0,-1);}\,}
\newkeycommand{\dkpointconj}[style=dkpointconj][1]{\,\tikz{\node[style=\commandkey{style}] (x) at (0,-0.1) {$#1$};\draw[boldedge] (x)--(0,1);}\,}

\newkeycommand{\blackkpoint}[style=black kpoint][1]{\,\tikz{\node[style=\commandkey{style}] (x) at (0,-0.1) {$#1$};\draw (x)--(0,1);}\,}
\newkeycommand{\blackkpointadj}[style=black kpointadj][1]{\,\tikz{\node[style=\commandkey{style}] (x) at (0,0.1) {$#1$};\draw (x)--(0,-1);}\,}
\newkeycommand{\blackdkpoint}[style=black dkpoint][1]{\,\tikz{\node[style=\commandkey{style}] (x) at (0,-0.1) {$#1$};\draw[boldedge] (x)--(0,1);}\,}
\newkeycommand{\blackdkpointadj}[style=black dkpointadj][1]{\,\tikz{\node[style=\commandkey{style}] (x) at (0,0.1) {$#1$};\draw[boldedge] (x)--(0,-1);}\,}

\newcommand{\trace}{\,\begin{tikzpicture}
  \begin{pgfonlayer}{nodelayer}
    \node [style=none] (0) at (0, -0.60) {};
    \node [style=upground] (1) at (0, 0.40) {};
  \end{pgfonlayer}
  \begin{pgfonlayer}{edgelayer}
    \draw [style=none] (0.center) to (1);
  \end{pgfonlayer}
\end{tikzpicture}\,}

\newcommand\discard\trace

\newkeycommand{\pointketbra}[style1=copoint,style2=point,estyle=][2]{\,%
\begin{tikzpicture}
  \begin{pgfonlayer}{nodelayer}
    \node [style=none] (0) at (0, -1.75) {};
    \node [style=\commandkey{style1}] (1) at (0, -0.75) {$#1$};
    \node [style=\commandkey{style2}] (2) at (0, 0.75) {$#2$};
    \node [style=none] (3) at (0, 1.75) {};
  \end{pgfonlayer}
  \begin{pgfonlayer}{edgelayer}
    \draw [\commandkey{estyle}] (0.center) to (1);
    \draw [\commandkey{estyle}] (2) to (3.center);
  \end{pgfonlayer}
\end{tikzpicture}\,}

\newkeycommand{\twocopointketbra}[style1=copoint,style2=copoint,style3=point][3]{\,%
\begin{tikzpicture}
  \begin{pgfonlayer}{nodelayer}
    \node [style=none] (0) at (-0.75, -1.75) {};
    \node [style=none] (1) at (0.75, -1.75) {};
    \node [style=\commandkey{style1}] (2) at (-0.75, -0.75) {$#1$};
    \node [style=\commandkey{style2}] (3) at (0.75, -0.75) {$#2$};
    \node [style=\commandkey{style3}] (4) at (0, 0.75) {$#3$};
    \node [style=none] (5) at (0, 1.75) {};
  \end{pgfonlayer}
  \begin{pgfonlayer}{edgelayer}
    \draw (0.center) to (2);
    \draw (1.center) to (3);
    \draw (4) to (5.center);
  \end{pgfonlayer}
\end{tikzpicture}\,}

\newkeycommand{\twopointketbra}[style1=copoint,style2=point,style3=point][3]{\,%
\begin{tikzpicture}
  \begin{pgfonlayer}{nodelayer}
    \node [style=none] (0) at (0, -1.75) {};
    \node [style=none] (1) at (-0.75, 1.75) {};
    \node [style=\commandkey{style1}] (2) at (0, -0.75) {$#1$};
    \node [style=\commandkey{style2}] (3) at (-0.75, 0.75) {$#2$};
    \node [style=\commandkey{style3}] (4) at (0.75, 0.75) {$#3$};
    \node [style=none] (5) at (0.75, 1.75) {};
  \end{pgfonlayer}
  \begin{pgfonlayer}{edgelayer}
    \draw (0.center) to (2);
    \draw (1.center) to (3);
    \draw (4) to (5.center);
  \end{pgfonlayer}
\end{tikzpicture}\,}

\newkeycommand{\pointbraket}[style1=point,style2=copoint,estyle=][2]{\,%
\begin{tikzpicture}
  \begin{pgfonlayer}{nodelayer}
    \node [style=\commandkey{style1}] (0) at (0, -0.625) {$#1$};
    \node [style=\commandkey{style2}] (1) at (0, 0.625) {$#2$};
  \end{pgfonlayer}
  \begin{pgfonlayer}{edgelayer}
    \draw [\commandkey{estyle}] (0) to (1);
  \end{pgfonlayer}
\end{tikzpicture}\,}

\newkeycommand{\kpointketbra}[style1=kpointdag,style2=kpoint,estyle=][2]{\,%
\begin{tikzpicture}
  \begin{pgfonlayer}{nodelayer}
    \node [style=none] (0) at (0, -1.9) {};
    \node [style=\commandkey{style1}] (1) at (0, -0.95) {$#1$};
    \node [style=\commandkey{style2}] (2) at (0, 0.95) {$#2$};
    \node [style=none] (3) at (0, 1.9) {};
  \end{pgfonlayer}
  \begin{pgfonlayer}{edgelayer}
    \draw [\commandkey{estyle}] (0.center) to (1);
    \draw [\commandkey{estyle}] (2) to (3.center);
  \end{pgfonlayer}
\end{tikzpicture}\,}

\newkeycommand{\kpointmap}[style1=kpoint,style2=map][2]{\,%
\begin{tikzpicture}
  \begin{pgfonlayer}{nodelayer}
    \node [style=\commandkey{style1}] (0) at (0, -1.1) {$#1$};
    \node [style=\commandkey{style2}] (1) at (0, 0.5) {$#2$};
    \node [style=none] (2) at (0, 1.75) {};
  \end{pgfonlayer}
  \begin{pgfonlayer}{edgelayer}
    \draw (0) to (1);
    \draw (1) to (2.center);
  \end{pgfonlayer}
\end{tikzpicture}%
\,}

\newkeycommand{\kpointmaptrans}[style1=kpointconj,style2=mapadj][2]{\,%
\begin{tikzpicture}
  \begin{pgfonlayer}{nodelayer}
    \node [style=\commandkey{style1}] (0) at (0, -1.1) {$#1$};
    \node [style=\commandkey{style2}] (1) at (0, 0.5) {$#2$};
    \node [style=none] (2) at (0, 1.75) {};
  \end{pgfonlayer}
  \begin{pgfonlayer}{edgelayer}
    \draw (0) to (1);
    \draw (1) to (2.center);
  \end{pgfonlayer}
\end{tikzpicture}%
\,}

\newkeycommand{\kpointmapadj}[style1=kpointadj,style2=map][2]{\,%
\begin{tikzpicture}
  \begin{pgfonlayer}{nodelayer}
    \node [style=\commandkey{style1}] (0) at (0, 1.1) {$#1$};
    \node [style=\commandkey{style2}] (1) at (0, -0.5) {$#2$};
    \node [style=none] (2) at (0, -1.75) {};
  \end{pgfonlayer}
  \begin{pgfonlayer}{edgelayer}
    \draw (0) to (1);
    \draw (1) to (2.center);
  \end{pgfonlayer}
\end{tikzpicture}%
\,}

\newkeycommand{\dkpointmap}[style1=dkpoint,style2=dmap][2]{\,%
\begin{tikzpicture}
  \begin{pgfonlayer}{nodelayer}
    \node [style=\commandkey{style1}] (0) at (0, -1.2) {$#1$};
    \node [style=\commandkey{style2}] (1) at (0, 0.4) {$#2$};
    \node [style=none] (2) at (0, 1.7) {};
  \end{pgfonlayer}
  \begin{pgfonlayer}{edgelayer}
    \draw[style=boldedge] (0) to (1);
    \draw[style=boldedge] (1) to (2.center);
  \end{pgfonlayer}
\end{tikzpicture}%
\,}

\newkeycommand{\dkpointmapx}[style1=dkpoint,style2=dmap][2]{\,%
\begin{tikzpicture}
  \begin{pgfonlayer}{nodelayer}
    \node [style=\commandkey{style1}] (0) at (0, -1.4) {$#1$};
    \node [style=\commandkey{style2}] (1) at (0, 0.6) {$#2$};
    \node [style=none] (2) at (0, 1.9) {};
  \end{pgfonlayer}
  \begin{pgfonlayer}{edgelayer}
    \draw[style=boldedge] (0) to (1);
    \draw[style=boldedge] (1) to (2.center);
  \end{pgfonlayer}
\end{tikzpicture}%
\,}

\newkeycommand{\boxcopointmapdag}[style1=map,style2=copoint][2]{\,%
\begin{tikzpicture}
  \begin{pgfonlayer}{nodelayer}
    \node [style=none] (0) at (0, -1.75) {};
    \node [style=\commandkey{style1}] (1) at (0, -0.5) {$#1$};
    \node [style=\commandkey{style2}] (2) at (0, 1) {$#2$};
  \end{pgfonlayer}
  \begin{pgfonlayer}{edgelayer}
    \draw (0.center) to (1);
    \draw (1) to (2);
  \end{pgfonlayer}
\end{tikzpicture}%
\,}

\newkeycommand{\kpointdagmap}[style1=map,style2=kpointdag][2]{\,%
\begin{tikzpicture}
  \begin{pgfonlayer}{nodelayer}
    \node [style=none] (0) at (0, -1.75) {};
    \node [style=\commandkey{style1}] (1) at (0, -0.5) {$#1$};
    \node [style=\commandkey{style2}] (2) at (0, 1.1) {$#2$};
  \end{pgfonlayer}
  \begin{pgfonlayer}{edgelayer}
    \draw (0.center) to (1);
    \draw (1) to (2);
  \end{pgfonlayer}
\end{tikzpicture}%
\,}

\newkeycommand{\sandwichmap}[style1=point,style2=map,style3=copoint][3]{\,%
\begin{tikzpicture}
  \begin{pgfonlayer}{nodelayer}
    \node [style=\commandkey{style1}] (0) at (0, -1.50) {$#1$};
    \node [style=\commandkey{style2}] (1) at (0, 0) {$#2$};
    \node [style=\commandkey{style3}] (2) at (0, 1.50) {$#3$};
  \end{pgfonlayer}
  \begin{pgfonlayer}{edgelayer}
    \draw (0) to (1);
    \draw (1) to (2);
  \end{pgfonlayer}
\end{tikzpicture}%
}

\newkeycommand{\kpointsandwichmap}[style1=kpoint,style2=map,style3=kpointdag][3]{\,%
\begin{tikzpicture}
  \begin{pgfonlayer}{nodelayer}
    \node [style=\commandkey{style1}] (0) at (0, -1.5) {$#1$};
    \node [style=\commandkey{style2}] (1) at (0, 0) {$#2$};
    \node [style=\commandkey{style3}] (2) at (0, 1.5) {$#3$};
  \end{pgfonlayer}
  \begin{pgfonlayer}{edgelayer}
    \draw (0) to (1);
    \draw (1) to (2);
  \end{pgfonlayer}
\end{tikzpicture}%
}

\newkeycommand{\sandwichmapconj}[style1=point,style2=mapconj,style3=copoint][3]{\,%
\begin{tikzpicture}
  \begin{pgfonlayer}{nodelayer}
    \node [style=\commandkey{style1}] (0) at (0, -1.50) {$#1$};
    \node [style=\commandkey{style2}] (1) at (0, 0) {$#2$};
    \node [style=\commandkey{style3}] (2) at (0, 1.50) {$#3$};
  \end{pgfonlayer}
  \begin{pgfonlayer}{edgelayer}
    \draw (0) to (1);
    \draw (1) to (2);
  \end{pgfonlayer}
\end{tikzpicture}%
}

\newkeycommand{\sandwichtwo}[style1=point,style2=map,style3=map,style4=copoint][4]{\,%
\begin{tikzpicture}
  \begin{pgfonlayer}{nodelayer}
    \node [style=\commandkey{style2}] (0) at (0, -1) {$#2$};
    \node [style=\commandkey{style3}] (1) at (0, 1) {$#3$};
    \node [style=\commandkey{style1}] (2) at (0, -2.6) {$#1$};
    \node [style=\commandkey{style4}] (3) at (0, 2.6) {$#4$};
  \end{pgfonlayer}
  \begin{pgfonlayer}{edgelayer}
    \draw (2) to (0);
    \draw (0) to (1);
    \draw (1) to (3);
  \end{pgfonlayer}
\end{tikzpicture}\,}

\newkeycommand{\longbraket}[style1=point,style2=copoint][2]{\,%
\begin{tikzpicture}
  \begin{pgfonlayer}{nodelayer}
    \node [style=\commandkey{style1}] (0) at (0, -2.6) {$#1$};
    \node [style=\commandkey{style2}] (1) at (0, 2.6) {$#2$};
  \end{pgfonlayer}
  \begin{pgfonlayer}{edgelayer}
    \draw (0) to (1);
  \end{pgfonlayer}
\end{tikzpicture}\,}

\newkeycommand{\onb}[style=point]{\ensuremath{\{\,\tikz{\node[style=\commandkey{style}] (x) at (0,-0.3) {$j$};\draw (x)--(0,0.7);}\,\}}\xspace}

\newkeycommand{\phase}[style=white phase dot][1]{\,\begin{tikzpicture}
    \begin{pgfonlayer}{nodelayer}
        \node [style=none] (0) at (0, 1) {};
        \node [style=\commandkey{style}] (2) at (0, -0) {$#1$}; 
        \node [style=none] (3) at (0, -1) {};
    \end{pgfonlayer}
    \begin{pgfonlayer}{edgelayer}
        \draw (2) to (0.center);
        \draw (3.center) to (2);
    \end{pgfonlayer}
\end{tikzpicture}\,}

\newkeycommand{\dphase}[style=white phase ddot][1]{\,\begin{tikzpicture}
    \begin{pgfonlayer}{nodelayer}
        \node [style=none] (0) at (0, 1.25) {};
        \node [style=\commandkey{style}] (2) at (0, -0) {$#1$};
        \node [style=none] (3) at (0, -1.25) {};
    \end{pgfonlayer}
    \begin{pgfonlayer}{edgelayer}
        \draw [boldedge] (2) to (0.center);
        \draw [boldedge] (3.center) to (2);
    \end{pgfonlayer}
\end{tikzpicture}\,}

\newkeycommand{\dphasepoint}[style=white phase ddot][1]{\,\begin{tikzpicture}[yshift=-3mm]
  \begin{pgfonlayer}{nodelayer}
    \node [style=none] (0) at (0, 1) {};
    \node [style=\commandkey{style}] (1) at (0, 0) {$#1$};
  \end{pgfonlayer}
  \begin{pgfonlayer}{edgelayer}
    \draw [style=boldedge] (0.center) to (1);
  \end{pgfonlayer}
\end{tikzpicture}\,}

\newkeycommand{\dphasepointgray}[style=gray phase ddot][1]{\,\begin{tikzpicture}[yshift=-3mm]
  \begin{pgfonlayer}{nodelayer}
    \node [style=none] (0) at (0, 1) {};
    \node [style=\commandkey{style}] (1) at (0, 0) {$#1$};
  \end{pgfonlayer}
  \begin{pgfonlayer}{edgelayer}
    \draw [style=boldedge] (0.center) to (1);
  \end{pgfonlayer}
\end{tikzpicture}\,}

\newkeycommand{\dphasecopoint}[style=white phase ddot][1]{\,\begin{tikzpicture}[yshift=3mm,yscale=-1]
  \begin{pgfonlayer}{nodelayer}
    \node [style=none] (0) at (0, 1) {};
    \node [style=\commandkey{style}] (1) at (0, 0) {$#1$};
  \end{pgfonlayer}
  \begin{pgfonlayer}{edgelayer}
    \draw [style=boldedge] (0.center) to (1);
  \end{pgfonlayer}
\end{tikzpicture}\,}

\newkeycommand{\dphasepointsm}[style=white phase ddot][1]{\,\begin{tikzpicture}[yshift=-3mm]
  \begin{pgfonlayer}{nodelayer}
    \node [style=none] (0) at (0, 1) {};
    \node [style=\commandkey{style}] (1) at (0, 0) {\footnotesize\!$#1$\!};
  \end{pgfonlayer}
  \begin{pgfonlayer}{edgelayer}
    \draw [style=boldedge] (0.center) to (1);
  \end{pgfonlayer}
\end{tikzpicture}\,}

\newkeycommand{\phasepoint}[style=white phase dot][1]{\,\begin{tikzpicture}[yshift=-3mm]
  \begin{pgfonlayer}{nodelayer}
    \node [style=none] (0) at (0, 1) {};
    \node [style=\commandkey{style}] (1) at (0, 0) {$#1$};
  \end{pgfonlayer}
  \begin{pgfonlayer}{edgelayer}
    \draw (0.center) to (1);
  \end{pgfonlayer}
\end{tikzpicture}\,}

\newkeycommand{\phasecopoint}[style=white phase dot][1]{\,\begin{tikzpicture}[yshift=3mm]
  \begin{pgfonlayer}{nodelayer}
    \node [style=none] (0) at (0, -1) {};
    \node [style=\commandkey{style}] (1) at (0, 0) {$#1$};
  \end{pgfonlayer}
  \begin{pgfonlayer}{edgelayer}
    \draw (0.center) to (1);
  \end{pgfonlayer}
\end{tikzpicture}\,}

\newkeycommand{\kpointbraket}[style1=kpoint,style2=kpoint adjoint,estyle=][2]{\,%
\begin{tikzpicture}
  \begin{pgfonlayer}{nodelayer}
    \node [style=\commandkey{style1}] (0) at (0, -0.735) {$#1$};
    \node [style=\commandkey{style2}] (1) at (0, 0.735) {$#2$};
  \end{pgfonlayer}
  \begin{pgfonlayer}{edgelayer}
    \draw [\commandkey{estyle}] (0) to (1);
  \end{pgfonlayer}
\end{tikzpicture}\,}

\newkeycommand{\kkpointbraket}[style1=kpoint,style2=kpoint adjoint,estyle=][2]{\,%
\begin{tikzpicture}
  \begin{pgfonlayer}{nodelayer}
    \node [style=\commandkey{style1}] (0) at (0, -0.685) {$#1$};
    \node [style=\commandkey{style2}] (1) at (0, 0.685) {$#2$};
  \end{pgfonlayer}
  \begin{pgfonlayer}{edgelayer}
    \draw [\commandkey{estyle}] (0) to (1);
  \end{pgfonlayer}
\end{tikzpicture}\,}

\newkeycommand{\kpointbraketx}[style1=kpoint,style2=kpoint adjoint,estyle=][2]{\,%
\begin{tikzpicture}
  \begin{pgfonlayer}{nodelayer}
    \node [style=\commandkey{style1}] (0) at (0, -0.885) {$#1$};
    \node [style=\commandkey{style2}] (1) at (0, 0.885) {$#2$};
  \end{pgfonlayer}
  \begin{pgfonlayer}{edgelayer}
    \draw [\commandkey{estyle}] (0) to (1);
  \end{pgfonlayer}
\end{tikzpicture}\,}

\newcommand{\bra}[1]{\ensuremath{\left\langle #1 \right|}}
\newcommand{\ket}[1]{\ensuremath{\left|  #1 \right\rangle}}

\tikzstyle{cdiag}=[matrix of math nodes, row sep=3em, column sep=3em, text height=1.5ex, text depth=0.25ex,inner sep=0.5em]
\tikzstyle{arrow above}=[transform canvas={yshift=0.5ex}]
\tikzstyle{arrow below}=[transform canvas={yshift=-0.5ex}]



%% file: macros/tikzstyles.tex
\usepackage{tikz}
\usetikzlibrary{decorations.markings}
\usetikzlibrary{shapes.geometric}

\pgfdeclarelayer{edgelayer}
\pgfdeclarelayer{nodelayer}
\pgfsetlayers{edgelayer,nodelayer,main}

\tikzstyle{none}=[inner sep=0pt]
\definecolor{hexcolor0xff0000}{rgb}{1.000,0.000,0.000}
\definecolor{hexcolor0x000000}{rgb}{0.000,0.000,0.000}
\definecolor{hexcolor0x00ff00}{rgb}{0.000,1.000,0.000}
\definecolor{hexcolor0x000000}{rgb}{0.000,0.000,0.000}
\definecolor{hexcolor0xffff00}{rgb}{1.000,1.000,0.000}
\definecolor{hexcolor0xffffff}{rgb}{1.000,1.000,1.000}

\tikzstyle{rn}=[circle,fill=hexcolor0xff0000,draw=hexcolor0x000000,line width=0.8 pt]
\tikzstyle{gn}=[circle,fill=hexcolor0x00ff00,draw=hexcolor0x000000,line width=0.8 pt]
\tikzstyle{yn}=[circle,fill=hexcolor0xffff00,draw=hexcolor0x000000,line width=0.8 pt]
\tikzstyle{wn}=[circle,fill=hexcolor0xffffff,draw=hexcolor0x000000,line width=0.8 pt]
\tikzstyle{wnthick}=[circle,fill=hexcolor0xffffff,draw=hexcolor0x000000,line width=2.500]

\tikzstyle{simple}=[-,draw=hexcolor0x000000,line width=2.000]
\tikzstyle{arrow}=[-,draw=hexcolor0x000000,postaction={decorate},decoration={markings,mark=at position .5 with {\arrow{>}}},line width=2.000]
\tikzstyle{tick}=[-,draw=hexcolor0x000000,postaction={decorate},decoration={markings,mark=at position .5 with {\draw (0,-0.1) -- (0,0.1);}},line width=2.000]
\tikzstyle{halfthickness}=[-,draw=hexcolor0x000000,line width=0.500]
\tikzstyle{thick}=[-,draw=hexcolor0x000000,line width=2.500]
\tikzstyle{thicker}=[-,draw=hexcolor0x000000,line width=4.000]

\tikzstyle{env}=[copoint,regular polygon rotate=0,minimum width=0.2cm, fill=black]

\tikzstyle{probs}=[shape=semicircle,fill=white,draw=black,shape border rotate=180,minimum width=1.2cm]

\tikzstyle{every picture}=[baseline=-0.25em,scale=0.5]
\tikzstyle{dotpic}=[] 
\tikzstyle{diredges}=[every to/.style={diredge}]
\tikzstyle{math matrix}=[matrix of math nodes,left delimiter=(,right delimiter=),inner sep=2pt,column sep=1em,row sep=0.5em,nodes={inner sep=0pt},text height=1.5ex, text depth=0.25ex]

\tikzstyle{inline text}=[text height=1.5ex, text depth=0.25ex,yshift=0.5mm]
\tikzstyle{label}=[font=\footnotesize,text height=1.5ex, text depth=0.25ex,yshift=0.5mm]
\tikzstyle{left label}=[label,anchor=east,xshift=1.5mm]
\tikzstyle{right label}=[label,anchor=west,xshift=-1.5mm]

\tikzstyle{braceedge}=[decorate,decoration={brace,amplitude=2mm,raise=-1mm}]
\tikzstyle{small braceedge}=[decorate,decoration={brace,amplitude=1mm,raise=-1mm}]

\tikzstyle{doubled}=[line width=1.6pt] 
\tikzstyle{boldedge}=[doubled,shorten <=-0.17mm,shorten >=-0.17mm]
\tikzstyle{boldedgegray}=[doubled,gray,shorten <=-0.17mm,shorten >=-0.17mm]

\tikzstyle{semidoubled}=[line width=1.4pt] 
\tikzstyle{semiboldedgegray}=[semidoubled,gray,shorten <=-0.17mm,shorten >=-0.17mm]

\tikzstyle{boldedgedashed}=[very thick,dashed,shorten <=-0.17mm,shorten >=-0.17mm]
\tikzstyle{vboldedgedashed}=[doubled,dashed,shorten <=-0.17mm,shorten >=-0.17mm]
\tikzstyle{left hook arrow}=[left hook-latex]
\tikzstyle{right hook arrow}=[right hook-latex]
\tikzstyle{sembracket}=[line width=0.5pt,shorten <=-0.07mm,shorten >=-0.07mm]

\tikzstyle{causal edge}=[->,thick,gray]
\tikzstyle{causal nondir}=[thick,gray]
\tikzstyle{timeline}=[thick,gray, dashed]

\tikzstyle{cedge}=[<->,thick,gray!70!white]

\tikzstyle{empty diagram}=[draw=gray!40!white,dashed,shape=rectangle,minimum width=1cm,minimum height=1cm]
\tikzstyle{empty diagram small}=[draw=gray!50!white,dashed,shape=rectangle,minimum width=0.6cm,minimum height=0.5cm]

\tikzstyle{dot}=[inner sep=0mm,minimum width=2mm,minimum height=2mm,draw,shape=circle]
\tikzstyle{ddot}=[inner sep=0mm, doubled, minimum width=2.5mm,minimum height=2.5mm,draw,shape=circle]

\tikzstyle{black dot}=[dot,fill=black]
\tikzstyle{white dot}=[dot,fill=white,,text depth=-0.2mm]
\tikzstyle{green dot}=[white dot] 
\tikzstyle{gray dot}=[dot,fill=gray!40!white,,text depth=-0.2mm]
\tikzstyle{red dot}=[gray dot] 

\tikzstyle{black ddot}=[ddot,fill=black]
\tikzstyle{white ddot}=[ddot,fill=white]
\tikzstyle{gray ddot}=[ddot,fill=gray!40!white]

\tikzstyle{gray edge}=[gray!40!white]

\tikzstyle{small dot}=[inner sep=0.5mm,minimum width=0pt,minimum height=0pt,draw,shape=circle]

\tikzstyle{small black dot}=[small dot,fill=black]
\tikzstyle{small white dot}=[small dot,fill=white]
\tikzstyle{small gray dot}=[small dot,fill=gray!40!white]

\tikzstyle{causal dot}=[inner sep=0.4mm,minimum width=0pt,minimum height=0pt,draw=white,shape=circle,fill=gray!40!white]

\tikzstyle{phase dimensions}=[minimum size=5mm,font=\footnotesize,rectangle,rounded corners=2.5mm,inner sep=0.2mm,outer sep=-2mm]
\tikzstyle{dphase dimensions}=[minimum size=5mm,font=\footnotesize,rectangle,rounded corners=2.5mm,inner sep=0.2mm,outer sep=-2mm]
\tikzstyle{phase dimensions small}=[minimum size=3.0mm,font=\footnotesize,rectangle,rounded corners=1.5mm,inner sep=0.2mm,outer sep=-1.2mm]

\tikzstyle{white phase dot}=[dot,fill=white,phase dimensions]
\tikzstyle{white phase ddot}=[ddot,fill=white,dphase dimensions]
\tikzstyle{green phase ddot}=[ddot,fill=green,dphase dimensions]
\tikzstyle{white phase dot small}=[dot,fill=white,phase dimensions small]
\tikzstyle{gray phase dot small}=[dot,fill=gray!40!white,phase dimensions small]

\tikzstyle{white rect ddot}=[draw=black,fill=white,doubled,minimum size=5mm,font=\footnotesize,rectangle,rounded corners=2.5mm,inner sep=0.2mm]
\tikzstyle{gray rect ddot}=[draw=black,fill=gray!40!white,doubled,minimum size=6mm,font=\footnotesize,rectangle,rounded corners=3mm]

\tikzstyle{gray phase dot}=[dot,fill=gray!40!white,phase dimensions]
\tikzstyle{gray phase ddot}=[ddot,fill=gray!40!white,dphase dimensions]
\tikzstyle{red phase ddot}=[ddot,fill=red,dphase dimensions]

\tikzstyle{grey phase dot}=[gray phase dot]
\tikzstyle{grey phase ddot}=[gray phase ddot]

\tikzstyle{small phase dimensions}=[minimum size=4mm,font=\tiny,rectangle,rounded corners=2mm,inner sep=0.2mm,outer sep=-2mm]
\tikzstyle{small dphase dimensions}=[minimum size=4mm,font=\tiny,rectangle,rounded corners=2mm,inner sep=0.2mm,outer sep=-2mm]

\tikzstyle{small gray phase dot}=[dot,fill=gray!40!white,small phase dimensions]
\tikzstyle{small gray phase ddot}=[ddot,fill=gray!40!white,small dphase dimensions]

\tikzstyle{small map}=[draw,shape=rectangle,minimum height=4mm,minimum width=4mm,fill=white]

\tikzstyle{cnot}=[fill=white,shape=circle,inner sep=-1.4pt]

\tikzstyle{asym hadamard}=[fill=white,draw,shape=NEbox,inner sep=0.6mm,font=\footnotesize,minimum height=4mm]
\tikzstyle{asym hadamard conj}=[fill=white,draw,shape=NWbox,inner sep=0.6mm,font=\footnotesize,minimum height=4mm]
\tikzstyle{asym hadamard dag}=[fill=white,draw,shape=SEbox,inner sep=0.6mm,font=\footnotesize,minimum height=4mm]

\tikzstyle{hadamard}=[fill=white,draw,inner sep=0.6mm,font=\footnotesize,minimum height=4mm,minimum width=4mm]
\tikzstyle{small hadamard}=[fill=white,draw,inner sep=0.6mm,minimum height=1.5mm,minimum width=1.5mm]
\tikzstyle{dhadamard}=[hadamard,doubled]
\tikzstyle{small dhadamard}=[small hadamard,doubled]
\tikzstyle{small dhadamard rotate}=[small hadamard,doubled,rotate=45]
\tikzstyle{antipode}=[white dot,inner sep=0.3mm,font=\footnotesize]

\tikzstyle{scalar}=[diamond,draw,inner sep=0.5pt,font=\small]
\tikzstyle{dscalar}=[diamond,doubled, draw,inner sep=0.5pt,font=\small]

\tikzstyle{small box}=[rectangle,inline text,fill=white,draw,minimum height=5mm,yshift=-0.5mm,minimum width=5mm,font=\small]
\tikzstyle{small gray box}=[small box,fill=gray!30]
\tikzstyle{medium box}=[rectangle,inline text,fill=white,draw,minimum height=5mm,yshift=-0.5mm,minimum width=10mm,font=\small]
\tikzstyle{square box}=[small box] 
\tikzstyle{medium gray box}=[small box,fill=gray!30]
\tikzstyle{semilarge box}=[rectangle,inline text,fill=white,draw,minimum height=5mm,yshift=-0.5mm,minimum width=12.5mm,font=\small]
\tikzstyle{large box}=[rectangle,inline text,fill=white,draw,minimum height=5mm,yshift=-0.5mm,minimum width=15mm,font=\small]
\tikzstyle{large gray box}=[small box,fill=gray!30]

\tikzstyle{Bayes box}=[rectangle,fill=black,draw, minimum height=3mm, minimum width=3mm]

\tikzstyle{gray square point}=[small box,fill=gray!50]

\tikzstyle{dphase box white}=[dhadamard]
\tikzstyle{dphase box gray}=[dhadamard,fill=gray!50!white]

\tikzstyle{point}=[regular polygon,regular polygon sides=3,draw,scale=0.75,inner sep=-0.5pt,minimum width=9mm,fill=white,regular polygon rotate=180]
\tikzstyle{copoint}=[regular polygon,regular polygon sides=3,draw,scale=0.75,inner sep=-0.5pt,minimum width=9mm,fill=white]
\tikzstyle{dpoint}=[point,doubled]
\tikzstyle{dcopoint}=[copoint,doubled]

\tikzstyle{wide copoint}=[fill=white,draw,shape=isosceles triangle,shape border rotate=90,isosceles triangle stretches=true,inner sep=0pt,minimum width=1.5cm,minimum height=6.12mm]
\tikzstyle{wide point}=[fill=white,draw,shape=isosceles triangle,shape border rotate=-90,isosceles triangle stretches=true,inner sep=0pt,minimum width=1.5cm,minimum height=6.12mm,yshift=-0.0mm]
\tikzstyle{wide point plus}=[fill=white,draw,shape=isosceles triangle,shape border rotate=-90,isosceles triangle stretches=true,inner sep=0pt,minimum width=1.74cm,minimum height=7mm,yshift=-0.0mm]

\tikzstyle{wide dpoint}=[fill=white,doubled,draw,shape=isosceles triangle,shape border rotate=-90,isosceles triangle stretches=true,inner sep=0pt,minimum width=1.5cm,minimum height=6.12mm,yshift=-0.0mm]
\tikzstyle{wide dcopoint}=[fill=white,doubled,draw,shape=isosceles triangle,shape border rotate=90,isosceles triangle stretches=true,inner sep=0pt,minimum width=1.5cm,minimum height=6.12mm,yshift=-0.0mm]

\tikzstyle{tinypoint}=[regular polygon,regular polygon sides=3,draw,scale=0.55,inner sep=-0.15pt,minimum width=6mm,fill=white,regular polygon rotate=180]

\tikzstyle{white point}=[point]
\tikzstyle{white dpoint}=[dpoint]
\tikzstyle{green point}=[white point] 
\tikzstyle{white copoint}=[copoint]
\tikzstyle{gray point}=[point,fill=gray!40!white]
\tikzstyle{gray dpoint}=[gray point,doubled]
\tikzstyle{red point}=[gray point] 
\tikzstyle{gray copoint}=[copoint,fill=gray!40!white]
\tikzstyle{gray dcopoint}=[gray copoint,doubled]

\tikzstyle{white point guide}=[regular polygon,regular polygon sides=3,font=\scriptsize,draw,scale=0.65,inner sep=-0.5pt,minimum width=9mm,fill=white,regular polygon rotate=180]

\tikzstyle{black point}=[point,fill=black,font=\color{white}]
\tikzstyle{black copoint}=[copoint,fill=black,font=\color{white}]

\tikzstyle{tiny gray point}=[tinypoint,fill=gray!40!white]

\tikzstyle{diredge}=[->]
\tikzstyle{ddiredge}=[<->]
\tikzstyle{rdiredge}=[<-]
\tikzstyle{thickdiredge}=[->, very thick]
\tikzstyle{pointer edge}=[->,very thick,gray]
\tikzstyle{pointer edge part}=[very thick,gray]
\tikzstyle{dashed edge}=[dashed]
\tikzstyle{thick dashed edge}=[very thick,dashed]
\tikzstyle{thick gray dashed edge}=[thick dashed edge,gray!40]
\tikzstyle{thick map edge}=[very thick,|->]

\makeatletter
\newcommand{\boxshape}[3]{%
\pgfdeclareshape{#1}{
\inheritsavedanchors[from=rectangle] 
\inheritanchorborder[from=rectangle]
\inheritanchor[from=rectangle]{center}
\inheritanchor[from=rectangle]{north}
\inheritanchor[from=rectangle]{south}
\inheritanchor[from=rectangle]{west}
\inheritanchor[from=rectangle]{east}
\backgroundpath{
\southwest \pgf@xa=\pgf@x \pgf@ya=\pgf@y
\northeast \pgf@xb=\pgf@x \pgf@yb=\pgf@y

\@tempdima=#2
\@tempdimb=#3

\pgfpathmoveto{\pgfpoint{\pgf@xa - 5pt + \@tempdima}{\pgf@ya}}
\pgfpathlineto{\pgfpoint{\pgf@xa - 5pt - \@tempdima}{\pgf@yb}}
\pgfpathlineto{\pgfpoint{\pgf@xb + 5pt + \@tempdimb}{\pgf@yb}}
\pgfpathlineto{\pgfpoint{\pgf@xb + 5pt - \@tempdimb}{\pgf@ya}}
\pgfpathlineto{\pgfpoint{\pgf@xa - 5pt + \@tempdima}{\pgf@ya}}
\pgfpathclose
}
}}

\boxshape{NEbox}{0pt}{5pt}
\boxshape{SEbox}{0pt}{-5pt}
\boxshape{NWbox}{5pt}{0pt}
\boxshape{SWbox}{-5pt}{0pt}
\boxshape{EBox}{-3pt}{3pt}
\boxshape{WBox}{3pt}{-3pt}
\makeatother

\tikzstyle{cloud}=[shape=cloud,draw,minimum width=1.5cm,minimum height=1.5cm]

\tikzstyle{map}=[draw,shape=NEbox,inner sep=2pt,minimum height=6mm,fill=white]
\tikzstyle{dashedmap}=[draw,dashed,shape=NEbox,inner sep=2pt,minimum height=6mm,fill=white]
\tikzstyle{mapdag}=[draw,shape=SEbox,inner sep=2pt,minimum height=6mm,fill=white]
\tikzstyle{mapadj}=[draw,shape=SEbox,inner sep=2pt,minimum height=6mm,fill=white]
\tikzstyle{maptrans}=[draw,shape=SWbox,inner sep=2pt,minimum height=6mm,fill=white]
\tikzstyle{mapconj}=[draw,shape=NWbox,inner sep=2pt,minimum height=6mm,fill=white]

\tikzstyle{langmap}=[draw,shape=NEbox,inner sep=2pt,minimum height=2.4mm,minimum width=3.2mm,fill=white]
\tikzstyle{langmaptrans}=[draw,shape=SWbox,inner sep=2pt,minimum height=2.4mm,minimum width=3.2mm,fill=white]

\tikzstyle{medium map}=[draw,shape=NEbox,inner sep=2pt,minimum height=6mm,fill=white,minimum width=7mm]
\tikzstyle{medium map dag}=[draw,shape=SEbox,inner sep=2pt,minimum height=6mm,fill=white,minimum width=7mm]
\tikzstyle{medium map adj}=[draw,shape=SEbox,inner sep=2pt,minimum height=6mm,fill=white,minimum width=7mm]
\tikzstyle{medium map trans}=[draw,shape=SWbox,inner sep=2pt,minimum height=6mm,fill=white,minimum width=7mm]
\tikzstyle{medium map conj}=[draw,shape=NWbox,inner sep=2pt,minimum height=6mm,fill=white,minimum width=7mm]
\tikzstyle{semilarge map}=[draw,shape=NEbox,inner sep=2pt,minimum height=6mm,fill=white,minimum width=9.5mm]
\tikzstyle{semilarge map trans}=[draw,shape=SWbox,inner sep=2pt,minimum height=6mm,fill=white,minimum width=9.5mm]
\tikzstyle{semilarge map adj}=[draw,shape=SEbox,inner sep=2pt,minimum height=6mm,fill=white,minimum width=9.5mm]
\tikzstyle{semilarge map dag}=[draw,shape=SEbox,inner sep=2pt,minimum height=6mm,fill=white,minimum width=9.5mm]
\tikzstyle{semilarge map conj}=[draw,shape=NWbox,inner sep=2pt,minimum height=6mm,fill=white,minimum width=9.5mm]
\tikzstyle{large map}=[draw,shape=NEbox,inner sep=2pt,minimum height=6mm,fill=white,minimum width=12mm]
\tikzstyle{large map conj}=[draw,shape=NWbox,inner sep=2pt,minimum height=6mm,fill=white,minimum width=12mm]
\tikzstyle{very large map}=[draw,shape=NEbox,inner sep=2pt,minimum height=6mm,fill=white,minimum width=17mm]

\tikzstyle{medium dmap}=[draw,doubled,shape=NEbox,inner sep=2pt,minimum height=6mm,fill=white,minimum width=7mm]
\tikzstyle{medium dmap dag}=[draw,doubled,shape=SEbox,inner sep=2pt,minimum height=6mm,fill=white,minimum width=7mm]
\tikzstyle{medium dmap adj}=[draw,doubled,shape=SEbox,inner sep=2pt,minimum height=6mm,fill=white,minimum width=7mm]
\tikzstyle{medium dmap trans}=[draw,doubled,shape=SWbox,inner sep=2pt,minimum height=6mm,fill=white,minimum width=7mm]
\tikzstyle{medium dmap conj}=[draw,doubled,shape=NWbox,inner sep=2pt,minimum height=6mm,fill=white,minimum width=7mm]
\tikzstyle{semilarge dmap}=[draw,doubled,shape=NEbox,inner sep=2pt,minimum height=6mm,fill=white,minimum width=9.5mm]
\tikzstyle{semilarge dmap trans}=[draw,doubled,shape=SWbox,inner sep=2pt,minimum height=6mm,fill=white,minimum width=9.5mm]
\tikzstyle{semilarge dmap adj}=[draw,doubled,shape=SEbox,inner sep=2pt,minimum height=6mm,fill=white,minimum width=9.5mm]
\tikzstyle{semilarge dmap dag}=[draw,doubled,shape=SEbox,inner sep=2pt,minimum height=6mm,fill=white,minimum width=9.5mm]
\tikzstyle{semilarge dmap conj}=[draw,doubled,shape=NWbox,inner sep=2pt,minimum height=6mm,fill=white,minimum width=9.5mm]
\tikzstyle{large dmap}=[draw,doubled,shape=NEbox,inner sep=2pt,minimum height=6mm,fill=white,minimum width=12mm]
\tikzstyle{large dmap conj}=[draw,doubled,shape=NWbox,inner sep=2pt,minimum height=6mm,fill=white,minimum width=12mm]
\tikzstyle{large dmap trans}=[draw,doubled,shape=SWbox,inner sep=2pt,minimum height=6mm,fill=white,minimum width=12mm]
\tikzstyle{large dmap adj}=[draw,doubled,shape=SEbox,inner sep=2pt,minimum height=6mm,fill=white,minimum width=12mm]
\tikzstyle{large dmap dag}=[draw,doubled,shape=SEbox,inner sep=2pt,minimum height=6mm,fill=white,minimum width=12mm]
\tikzstyle{very large dmap}=[draw,doubled,shape=NEbox,inner sep=2pt,minimum height=6mm,fill=white,minimum width=19.5mm]

\tikzstyle{muxbox}=[draw,shape=rectangle,minimum height=3mm,minimum width=3mm,fill=white]
\tikzstyle{dmuxbox}=[muxbox,doubled]

\tikzstyle{box}=[draw,shape=rectangle,inner sep=2pt,minimum height=6mm,minimum width=6mm,fill=white]
\tikzstyle{dbox}=[draw,doubled,shape=rectangle,inner sep=2pt,minimum height=6mm,minimum width=6mm,fill=white]
\tikzstyle{dmap}=[draw,doubled,shape=NEbox,inner sep=2pt,minimum height=6mm,fill=white]
\tikzstyle{dmapdag}=[draw,doubled,shape=SEbox,inner sep=2pt,minimum height=6mm,fill=white]
\tikzstyle{dmapadj}=[draw,doubled,shape=SEbox,inner sep=2pt,minimum height=6mm,fill=white]
\tikzstyle{dmaptrans}=[draw,doubled,shape=SWbox,inner sep=2pt,minimum height=6mm,fill=white]
\tikzstyle{dmapconj}=[draw,doubled,shape=NWbox,inner sep=2pt,minimum height=6mm,fill=white]

\tikzstyle{ddmap}=[draw,doubled,dashed,shape=NEbox,inner sep=2pt,minimum height=6mm,fill=white]
\tikzstyle{ddmapdag}=[draw,doubled,dashed,shape=SEbox,inner sep=2pt,minimum height=6mm,fill=white]
\tikzstyle{ddmapadj}=[draw,doubled,dashed,shape=SEbox,inner sep=2pt,minimum height=6mm,fill=white]
\tikzstyle{ddmaptrans}=[draw,doubled,dashed,shape=SWbox,inner sep=2pt,minimum height=6mm,fill=white]
\tikzstyle{ddmapconj}=[draw,doubled,dashed,shape=NWbox,inner sep=2pt,minimum height=6mm,fill=white]

\boxshape{sNEbox}{0pt}{3pt}
\boxshape{sSEbox}{0pt}{-3pt}
\boxshape{sNWbox}{3pt}{0pt}
\boxshape{sSWbox}{-3pt}{0pt}
\tikzstyle{smap}=[draw,shape=sNEbox,fill=white]
\tikzstyle{smapdag}=[draw,shape=sSEbox,fill=white]
\tikzstyle{smapadj}=[draw,shape=sSEbox,fill=white]
\tikzstyle{smaptrans}=[draw,shape=sSWbox,fill=white]
\tikzstyle{smapconj}=[draw,shape=sNWbox,fill=white]

\tikzstyle{dsmap}=[draw,dashed,shape=sNEbox,fill=white]
\tikzstyle{dsmapdag}=[draw,dashed,shape=sSEbox,fill=white]
\tikzstyle{dsmaptrans}=[draw,dashed,shape=sSWbox,fill=white]
\tikzstyle{dsmapconj}=[draw,dashed,shape=sNWbox,fill=white]

\boxshape{mNEbox}{0pt}{10pt}
\boxshape{mSEbox}{0pt}{-10pt}
\boxshape{mNWbox}{10pt}{0pt}
\boxshape{mSWbox}{-10pt}{0pt}
\tikzstyle{mmap}=[draw,shape=mNEbox]
\tikzstyle{mmapdag}=[draw,shape=mSEbox]
\tikzstyle{mmaptrans}=[draw,shape=mSWbox]
\tikzstyle{mmapconj}=[draw,shape=mNWbox]

\tikzstyle{mmapgray}=[draw,fill=gray!40!white,shape=mNEbox]
\tikzstyle{smapgray}=[draw,fill=gray!40!white,shape=sNEbox]

\makeatletter
\pgfdeclareshape{cornerpoint}{
\inheritsavedanchors[from=rectangle] 
\inheritanchorborder[from=rectangle]
\inheritanchor[from=rectangle]{center}
\inheritanchor[from=rectangle]{north}
\inheritanchor[from=rectangle]{south}
\inheritanchor[from=rectangle]{west}
\inheritanchor[from=rectangle]{east}
\backgroundpath{
\southwest \pgf@xa=\pgf@x \pgf@ya=\pgf@y
\northeast \pgf@xb=\pgf@x \pgf@yb=\pgf@y

\pgfmathsetmacro{\pgf@shorten@left}{\pgfkeysvalueof{/tikz/shorten left}}
\pgfmathsetmacro{\pgf@shorten@right}{\pgfkeysvalueof{/tikz/shorten right}}

\pgfpathmoveto{\pgfpoint{0.5 * (\pgf@xa + \pgf@xb)}{\pgf@ya - 5pt}}
\pgfpathlineto{\pgfpoint{\pgf@xa - 8pt + \pgf@shorten@left}{\pgf@yb - 1.5 * \pgf@shorten@left}}
\pgfpathlineto{\pgfpoint{\pgf@xa - 8pt + \pgf@shorten@left}{\pgf@yb}}
\pgfpathlineto{\pgfpoint{\pgf@xb + 8pt - \pgf@shorten@right}{\pgf@yb}}
\pgfpathlineto{\pgfpoint{\pgf@xb + 8pt - \pgf@shorten@right}{\pgf@yb - 1.5 * \pgf@shorten@right}}
\pgfpathclose
}
}

\pgfdeclareshape{cornercopoint}{
\inheritsavedanchors[from=rectangle] 
\inheritanchorborder[from=rectangle]
\inheritanchor[from=rectangle]{center}
\inheritanchor[from=rectangle]{north}
\inheritanchor[from=rectangle]{south}
\inheritanchor[from=rectangle]{west}
\inheritanchor[from=rectangle]{east}
\backgroundpath{
\southwest \pgf@xa=\pgf@x \pgf@ya=\pgf@y
\northeast \pgf@xb=\pgf@x \pgf@yb=\pgf@y

\pgfmathsetmacro{\pgf@shorten@left}{\pgfkeysvalueof{/tikz/shorten left}}
\pgfmathsetmacro{\pgf@shorten@right}{\pgfkeysvalueof{/tikz/shorten right}}

\pgfpathmoveto{\pgfpoint{0.5 * (\pgf@xa + \pgf@xb)}{\pgf@yb + 5pt}}
\pgfpathlineto{\pgfpoint{\pgf@xa - 8pt + \pgf@shorten@left}{\pgf@ya + 1.5 * \pgf@shorten@left}}
\pgfpathlineto{\pgfpoint{\pgf@xa - 8pt + \pgf@shorten@left}{\pgf@ya}}
\pgfpathlineto{\pgfpoint{\pgf@xb + 8pt - \pgf@shorten@right}{\pgf@ya}}
\pgfpathlineto{\pgfpoint{\pgf@xb + 8pt - \pgf@shorten@right}{\pgf@ya + 1.5 * \pgf@shorten@right}}
\pgfpathclose
}
}

\pgfdeclareshape{langpoint}{
\inheritsavedanchors[from=rectangle] 
\inheritanchorborder[from=rectangle]
\inheritanchor[from=rectangle]{center}
\inheritanchor[from=rectangle]{north}
\inheritanchor[from=rectangle]{south}
\inheritanchor[from=rectangle]{west}
\inheritanchor[from=rectangle]{east}
\backgroundpath{
\southwest \pgf@xa=\pgf@x \pgf@ya=\pgf@y
\northeast \pgf@xb=\pgf@x \pgf@yb=\pgf@y

\pgfmathsetmacro{\pgf@shorten@left}{\pgfkeysvalueof{/tikz/shorten left}}
\pgfmathsetmacro{\pgf@shorten@right}{\pgfkeysvalueof{/tikz/shorten right}}

\pgfpathmoveto{\pgfpoint{0.5 * (\pgf@xa + \pgf@xb)}{\pgf@ya - 2pt}}
\pgfpathlineto{\pgfpoint{\pgf@xa - 8pt}{\pgf@yb - 3 * \pgf@shorten@left + 5pt}} 
\pgfpathlineto{\pgfpoint{\pgf@xa - 8pt}{\pgf@yb -1pt}}
\pgfpathlineto{\pgfpoint{\pgf@xb + 8pt}{\pgf@yb -1pt}}
\pgfpathlineto{\pgfpoint{\pgf@xb + 8pt}{\pgf@yb - 3 * \pgf@shorten@left + 5pt}}
\pgfpathclose
}
}

\pgfdeclareshape{langcopoint}{
\inheritsavedanchors[from=rectangle] 
\inheritanchorborder[from=rectangle]
\inheritanchor[from=rectangle]{center}
\inheritanchor[from=rectangle]{north}
\inheritanchor[from=rectangle]{south}
\inheritanchor[from=rectangle]{west}
\inheritanchor[from=rectangle]{east}
\backgroundpath{
\southwest \pgf@xa=\pgf@x \pgf@ya=\pgf@y
\northeast \pgf@xb=\pgf@x \pgf@yb=\pgf@y

\pgfmathsetmacro{\pgf@shorten@left}{\pgfkeysvalueof{/tikz/shorten left}}
\pgfmathsetmacro{\pgf@shorten@right}{\pgfkeysvalueof{/tikz/shorten right}}

\pgfpathmoveto{\pgfpoint{0.5 * (\pgf@xa + \pgf@xb)}{\pgf@yb +0pt}}
\pgfpathlineto{\pgfpoint{\pgf@xa - 8pt}{\pgf@ya + 3 * \pgf@shorten@left - 5pt}} 
\pgfpathlineto{\pgfpoint{\pgf@xa - 8pt}{\pgf@ya + 1pt}}
\pgfpathlineto{\pgfpoint{\pgf@xb + 8pt}{\pgf@ya + 1pt}}
\pgfpathlineto{\pgfpoint{\pgf@xb + 8pt}{\pgf@ya + 3 * \pgf@shorten@left - 5pt}}
\pgfpathclose
}
}

\pgfdeclareshape{langrect}{
\inheritsavedanchors[from=rectangle] 
\inheritanchorborder[from=rectangle]
\inheritanchor[from=rectangle]{center}
\inheritanchor[from=rectangle]{north}
\inheritanchor[from=rectangle]{south}
\inheritanchor[from=rectangle]{west}
\inheritanchor[from=rectangle]{east}
\backgroundpath{
\southwest \pgf@xa=\pgf@x \pgf@ya=\pgf@y
\northeast \pgf@xb=\pgf@x \pgf@yb=\pgf@y

\pgfmathsetmacro{\pgf@shorten@left}{\pgfkeysvalueof{/tikz/shorten left}}
\pgfmathsetmacro{\pgf@shorten@right}{\pgfkeysvalueof{/tikz/shorten right}}

\pgfpathmoveto{\pgfpoint{\pgf@xa - 8pt}{\pgf@ya + 3 * \pgf@shorten@left - 5pt}} 
\pgfpathlineto{\pgfpoint{\pgf@xa - 8pt}{\pgf@ya + 1pt}}
\pgfpathlineto{\pgfpoint{\pgf@xb + 8pt}{\pgf@ya + 1pt}}
\pgfpathlineto{\pgfpoint{\pgf@xb + 8pt}{\pgf@ya + 3 * \pgf@shorten@left - 5pt}}
\pgfpathclose
}
}

\makeatother

\pgfkeyssetvalue{/tikz/shorten left}{0pt}
\pgfkeyssetvalue{/tikz/shorten right}{0pt}

\tikzstyle{kpoint common}=[draw,fill=white,inner sep=1pt,minimum height=4mm]

\tikzstyle{langstate}=[shape=langcopoint,shorten left=5pt,kpoint common,font=\footnotesize]
\tikzstyle{langeffect}=[shape=langpoint,shorten left=5pt,kpoint common,font=\footnotesize]
\tikzstyle{langstatedash}=[shape=langcopoint,dashed, shorten left=5pt,kpoint common,font=\footnotesize]
\tikzstyle{langeffectdash}=[shape=langpoint,dashed, shorten left=5pt,kpoint common,font=\footnotesize]
\tikzstyle{langbox}=[shape=langrect,shorten left=5pt,kpoint common,font=\footnotesize] 

\tikzstyle{kpoint}=[shape=cornerpoint,shorten left=5pt,kpoint common]
\tikzstyle{kpoint adjoint}=[shape=cornercopoint,shorten left=5pt,kpoint common]

\tikzstyle{kpoint conjugate}=[shape=cornerpoint,shorten right=5pt,kpoint common]
\tikzstyle{kpoint transpose}=[shape=cornercopoint,shorten right=5pt,kpoint common]
\tikzstyle{kpoint symm}=[shape=cornerpoint,shorten left=5pt,shorten right=5pt,kpoint common]

\tikzstyle{black kpoint}=[shape=cornerpoint,shorten left=5pt,kpoint common,fill=black,font=\color{white}]
\tikzstyle{black kpoint adjoint}=[shape=cornercopoint,shorten left=5pt,kpoint common,fill=black,font=\color{white}]
\tikzstyle{black kpointadj}=[shape=cornercopoint,shorten left=5pt,kpoint common,fill=black,font=\color{white}]

\tikzstyle{black dkpoint}=[shape=cornerpoint,shorten left=5pt,kpoint common,fill=black, doubled,font=\color{white}]
\tikzstyle{black dkpoint adjoint}=[shape=cornercopoint,shorten left=5pt,kpoint common,fill=black, doubled,font=\color{white}]
\tikzstyle{black dkpointadj}=[shape=cornercopoint,shorten left=5pt,kpoint common,fill=black, doubled,font=\color{white}]

\tikzstyle{kpointdag}=[kpoint adjoint]
\tikzstyle{kpointadj}=[kpoint adjoint]
\tikzstyle{kpointconj}=[kpoint conjugate]
\tikzstyle{kpointtrans}=[kpoint transpose]

\tikzstyle{big kpoint}=[kpoint, minimum width=1.2 cm, minimum height=8mm, inner sep=4pt, text depth=3mm]

\tikzstyle{wide kpoint}=[kpoint, minimum width=1 cm, inner sep=2pt]
\tikzstyle{wide kpointdag}=[kpointdag, minimum width=1 cm, inner sep=2pt]
\tikzstyle{wide kpointconj}=[kpointconj, minimum width=1 cm, inner sep=2pt]
\tikzstyle{wide kpointtrans}=[kpointtrans, minimum width=1 cm, inner sep=2pt]

\tikzstyle{gray kpoint}=[kpoint,fill=gray!50!white]
\tikzstyle{gray kpointdag}=[kpointdag,fill=gray!50!white]
\tikzstyle{gray kpointadj}=[kpointadj,fill=gray!50!white]
\tikzstyle{gray kpointconj}=[kpointconj,fill=gray!50!white]
\tikzstyle{gray kpointtrans}=[kpointtrans,fill=gray!50!white]

\tikzstyle{gray dkpoint}=[kpoint,fill=gray!50!white,doubled]
\tikzstyle{gray dkpointdag}=[kpointdag,fill=gray!50!white,doubled]
\tikzstyle{gray dkpointadj}=[kpointadj,fill=gray!50!white,doubled]
\tikzstyle{gray dkpointconj}=[kpointconj,fill=gray!50!white,doubled]
\tikzstyle{gray dkpointtrans}=[kpointtrans,fill=gray!50!white,doubled]

\tikzstyle{white label}=[draw,fill=white,rectangle,inner sep=0.7 mm]
\tikzstyle{gray label}=[draw,fill=gray!50!white,rectangle,inner sep=0.7 mm]
\tikzstyle{black label}=[draw,fill=black,rectangle,inner sep=0.7 mm]

\tikzstyle{dkpoint}=[kpoint,doubled]
\tikzstyle{wide dkpoint}=[wide kpoint,doubled]
\tikzstyle{dkpointdag}=[kpoint adjoint,doubled]
\tikzstyle{wide dkpointdag}=[wide kpointdag,doubled]
\tikzstyle{dkcopoint}=[kpoint adjoint,doubled]
\tikzstyle{dkpointadj}=[kpoint adjoint,doubled]
\tikzstyle{dkpointconj}=[kpoint conjugate,doubled]
\tikzstyle{dkpointtrans}=[kpoint transpose,doubled]

\tikzstyle{kscalar}=[kpoint common, shape=EBox, inner xsep=-1pt, inner ysep=3pt,font=\small]
\tikzstyle{kscalarconj}=[kpoint common, shape=WBox, inner xsep=-1pt, inner ysep=3pt,font=\small]

 \tikzstyle{upground}=[circuit ee IEC,ground,rotate=90,scale=2.5]
 \tikzstyle{downground}=[circuit ee IEC,ground,rotate=-90,scale=2.5]
 \tikzstyle{bigground}=[regular polygon,regular polygon sides=3,draw=gray,scale=0.50,inner sep=-0.5pt,minimum width=10mm,fill=gray]

\tikzstyle{arrs}=[-latex,font=\small,auto]
\tikzstyle{arrow plain}=[arrs]
\tikzstyle{arrow dashed}=[dashed,arrs]
\tikzstyle{arrow bold}=[very thick,arrs]
\tikzstyle{arrow hide}=[draw=white!0,-]
\tikzstyle{arrow reverse}=[latex-]
\tikzstyle{cdnode}=[]

%% file: figures/Fan3.tikz
\begin{tikzpicture}[scale=0.80]
	\begin{pgfonlayer}{nodelayer}
		\node [style=none] (0) at (-2.5, -2) {};
		\node [style=white dot] (1) at (-3.5, -1) {};
		\node [style=none] (2) at (-4.5, -2) {};
		\node [style=none] (3) at (-3.5, -0.5) {};
		\node [style=none] (4) at (3, -2) {};
		\node [style=none] (5) at (-5, -2) {};
		\node [style=none] (6) at (3, 2.25) {};
		\node [style=none] (7) at (-5, 2.25) {};
		\node [style=none] (8) at (-1, 3.5) {};
		\node [style=none] (9) at (-3.5, -0.5) {};
		\node [style=none] (10) at (-1.5, -0.5) {};
		\node [style=none] (11) at (-0.75, 2.25) {};
		\node [style=white dot] (12) at (-2.5, 0.5) {};
		\node [style=none] (13) at (-2.5, 1.5) {};
		\node [style=none] (14) at (-4.25, -0.5) {};
		\node [style=none] (15) at (-2.5, 2.75) {};
		\node [style=none] (16) at (-3.5, -0.5) {};
		\node [style=langstate] (17) at (-2.5, 1.75) {\!\!\!\!{\tt black}\!\!\!\!};
		\node [style=none] (18) at (-0.75, -0.5) {};
		\node [style=none] (19) at (-4.25, 2.25) {};
		\node [style=none] (20) at (-1.5, -0.5) {};
		\node [style=none] (21) at (1.25, -0.5) {};
		\node [style=none] (22) at (1.25, -0.5) {};
		\node [style=langstate] (23) at (1.25, 1.75) {\!\!\!\!{\tt metal}\!\!\!\!};
		\node [style=none] (24) at (-4.5, -3.25) {};
		\node [style=none] (25) at (-2.5, -2) {};
		\node [style=langstate] (26) at (4.75, 1.75) {\!\!\!{\tt fan}\!\!\!};
		\node [style=none] (27) at (4.75, -2) {};
		\node [style=none] (28) at (4.75, -2) {};
	\end{pgfonlayer}
	\begin{pgfonlayer}{edgelayer}
		\draw [style=swap, in=-165, out=90, looseness=1.00] (2.center) to (1);
		\draw [style=swap, in=90, out=-15, looseness=1.00] (1) to (0.center);
		\draw [style=swap, dashed] (7.center) to (5.center);
		\draw [style=swap, dashed] (5.center) to (4.center);
		\draw [style=swap, dashed] (4.center) to (6.center);
		\draw [style=swap, dashed] (8.center) to (7.center);
		\draw [style=swap, dashed] (8.center) to (6.center);
		\draw [style=swap] (9.center) to (1);
		\draw [style=swap, in=-165, out=90, looseness=1.00] (16.center) to (12);
		\draw [style=swap, in=90, out=-15, looseness=1.00] (12) to (10.center);
		\draw [style=swap, dashed] (19.center) to (14.center);
		\draw [style=swap, dashed] (14.center) to (18.center);
		\draw [style=swap, dashed] (18.center) to (11.center);
		\draw [style=swap, dashed] (15.center) to (19.center);
		\draw [style=swap, dashed] (15.center) to (11.center);
		\draw [style=swap] (13.center) to (12);
		\draw [style=swap, bend right=90, looseness=1.00] (20.center) to (21.center);
		\draw [style=swap] (23) to (22.center);
		\draw [style=swap] (2.center) to (24.center);
		\draw [style=swap, bend right=90, looseness=0.50] (25.center) to (27.center);
		\draw [style=swap] (26) to (28.center);
	\end{pgfonlayer}
\end{tikzpicture}

%% file: figures/Fan2.tikz
\begin{tikzpicture}[scale=0.80]
	\begin{pgfonlayer}{nodelayer}
		\node [style=none] (0) at (-4.25, -2) {};
		\node [style=langstate] (1) at (-3.25, 1.25) {\!\!\!\!{\tt black}\!\!\!\!};
		\node [style=none] (2) at (-2.25, -1) {};
		\node [style=white dot] (3) at (-3.25, 0) {};
		\node [style=none] (4) at (-4.25, -1) {};
		\node [style=none] (5) at (0, -1) {};
		\node [style=none] (6) at (-1.5, -1) {};
		\node [style=none] (7) at (-5, -1) {};
		\node [style=none] (8) at (-1.5, 1.75) {};
		\node [style=none] (9) at (-5, 1.75) {};
		\node [style=none] (10) at (-3.25, 2.25) {};
		\node [style=none] (11) at (-3.25, 1) {};
		\node [style=none] (12) at (2, -1) {};
		\node [style=none] (13) at (2.75, 1.75) {};
		\node [style=white dot] (14) at (1, 0) {};
		\node [style=none] (15) at (1, 1) {};
		\node [style=none] (16) at (-0.75, -1) {};
		\node [style=none] (17) at (1, 2.25) {};
		\node [style=none] (18) at (0, -1) {};
		\node [style=langstate] (19) at (1, 1.25) {\!\!\!\!{\tt metal}\!\!\!\!};
		\node [style=none] (20) at (2.75, -1) {};
		\node [style=none] (21) at (-0.75, 1.75) {};
		\node [style=none] (22) at (2, -1) {};
		\node [style=none] (23) at (4.5, -1) {};
		\node [style=none] (24) at (4.5, -1) {};
		\node [style=langstate] (25) at (4.5, 1.25) {\!\!\!{\tt fan}\!\!\!};
	\end{pgfonlayer}
	\begin{pgfonlayer}{edgelayer}
		\draw [style=swap, in=-165, out=90, looseness=1.00] (4.center) to (3);
		\draw [style=swap, in=90, out=-15, looseness=1.00] (3) to (2.center);
		\draw [style=swap, bend right=90, looseness=1.25] (2.center) to (5.center);
		\draw [style=swap] (4.center) to (0.center);
		\draw [style=swap, dashed] (9.center) to (7.center);
		\draw [style=swap, dashed] (7.center) to (6.center);
		\draw [style=swap, dashed] (6.center) to (8.center);
		\draw [style=swap, dashed] (10.center) to (9.center);
		\draw [style=swap, dashed] (10.center) to (8.center);
		\draw [style=swap] (11.center) to (3);
		\draw [style=swap, in=-165, out=90, looseness=1.00] (18.center) to (14);
		\draw [style=swap, in=90, out=-15, looseness=1.00] (14) to (12.center);
		\draw [style=swap, dashed] (21.center) to (16.center);
		\draw [style=swap, dashed] (16.center) to (20.center);
		\draw [style=swap, dashed] (20.center) to (13.center);
		\draw [style=swap, dashed] (17.center) to (21.center);
		\draw [style=swap, dashed] (17.center) to (13.center);
		\draw [style=swap] (15.center) to (14);
		\draw [style=swap, bend right=90, looseness=1.00] (22.center) to (23.center);
		\draw [style=swap] (25) to (24.center);
	\end{pgfonlayer}
\end{tikzpicture}

%% file: figures/separate.tikz
\begin{tikzpicture}[scale=0.80]
	\begin{pgfonlayer}{nodelayer}
		\node [style=none] (0) at (-3.5, -0.75) {};
		\node [style=none] (1) at (-2, 0.25) {};
		\node [style=none] (2) at (4, 0.25) {};
		\node [style=none] (3) at (2, 0.25) {};
		\node [style=none] (4) at (-2, -0.75) {};
		\node [style=none] (5) at (-3.5, 0.25) {};
		\node [style=langstate] (6) at (2, 0.5) {\!\!\!\!$\psi_{A}$\!\!\!\!};
		\node [style=none] (7) at (4, -0.75) {};
		\node [style=none] (8) at (2, -0.75) {};
		\node [style=none] (9) at (0, 0) {$\not=$}; 
		\node [style=langstate] (10) at (4, 0.5) {\!\!\!\!$\psi_{B}$\!\!\!\!};
		\node [style=langstate] (11) at (-2.75, 0.5) {\!$\psi_{AB}$\!};
	\end{pgfonlayer}
	\begin{pgfonlayer}{edgelayer}
		\draw [style=swap] (4.center) to (1.center);
		\draw [style=swap] (0.center) to (5.center);
		\draw [style=swap] (8.center) to (3.center);
		\draw [style=swap] (7.center) to (2.center);
	\end{pgfonlayer}
\end{tikzpicture}

%% file: figures/box2.tikz
\begin{tikzpicture}[scale=0.80]
	\begin{pgfonlayer}{nodelayer}
		\node [style=none] (0) at (-1, -0.25) {};
		\node [style=none] (1) at (1, -1.25) {};
		\node [style=langbox] (2) at (0, 0) {\!\!\!{\tt process}\!\!\!};
		\node [style=none] (3) at (1, -0.25) {};
		\node [style=none] (4) at (-1, -1.25) {};
		\node [style=none] (5) at (-1, 0.25) {};
		\node [style=none] (6) at (-1, 1.25) {};
		\node [style=none] (7) at (1, 0.25) {};
		\node [style=none] (8) at (1, 1.25) {};
		\node [style=none] (9) at (0, 1) {\ldots};
		\node [style=none] (10) at (0, -1) {\ldots};
	\end{pgfonlayer}
	\begin{pgfonlayer}{edgelayer}
		\draw [style=swap] (1.center) to (3.center); 
		\draw [style=swap] (4.center) to (0.center);
		\draw [style=swap] (5.center) to (6.center);
		\draw [style=swap] (7.center) to (8.center);
	\end{pgfonlayer}
\end{tikzpicture}

%% file: figures/boxexample2.tikz
\begin{tikzpicture}[scale=0.80]
	\begin{pgfonlayer}{nodelayer}
		\node [style=none] (0) at (-1, -1) {};
		\node [style=none] (1) at (1.5, -3) {};
		\node [style=langbox] (2) at (0.25, -0.75) {\!\!\!{\tt process 2}\!\!\!};
		\node [style=none] (3) at (1.5, -1) {};
		\node [style=none] (4) at (-1, -3) {};
		\node [style=none] (5) at (-1, -0.5) {};
		\node [style=none] (6) at (-1, 0.5) {};
		\node [style=none] (7) at (1.5, -0.5) {};
		\node [style=none] (8) at (1.5, 0.5) {};
		\node [style=langbox] (9) at (2.75, 0.75) {\!\!\!{\tt process 1}\!\!\!};   
		\node [style=none] (10) at (4, -1.75) {};
		\node [style=none] (11) at (4, 0.5) {};
		\node [style=langstate] (12) at (-2, 0.75) {\!{\tt state}\!};
		\node [style=none] (13) at (-3, 0.5) {};
		\node [style=none] (14) at (-3, -3) {};
		\node [style=none] (15) at (1.5, 1) {};
		\node [style=none] (16) at (1.5, 1.75) {};
		\node [style=none] (17) at (4, 1) {};
		\node [style=none] (18) at (4, 1.75) {};
		\node [style=none] (19) at (6, 1.75) {};
		\node [style=none] (20) at (6, -1.75) {};
		\node [style=langeffect] (21) at (5, -2) {\!\!{\tt effect}\!\!};
	\end{pgfonlayer}
	\begin{pgfonlayer}{edgelayer}
		\draw [style=swap] (1.center) to (3.center);
		\draw [style=swap] (4.center) to (0.center);
		\draw [style=swap] (5.center) to (6.center);
		\draw [style=swap] (7.center) to (8.center);
		\draw [style=swap] (10.center) to (11.center);
		\draw [style=swap] (14.center) to (13.center);
		\draw [style=swap] (15.center) to (16.center);
		\draw [style=swap] (17.center) to (18.center);
		\draw [style=swap] (20.center) to (19.center);
	\end{pgfonlayer}
\end{tikzpicture}

%% file: figures/cap.tikz
\begin{tikzpicture}[scale=0.80]
	\begin{pgfonlayer}{nodelayer}
		\node [style=none] (0) at (1.25, -0.5) {};
		\node [style=none] (1) at (-1.25, -0.5) {};
	\end{pgfonlayer}
	\begin{pgfonlayer}{edgelayer}
		\draw [style=swap, bend left=90, looseness=1.50] (1.center) to (0.center);  
	\end{pgfonlayer} 
\end{tikzpicture}

%% file: figures/cup.tikz
\begin{tikzpicture}[scale=0.80]
	\begin{pgfonlayer}{nodelayer}
		\node [style=none] (0) at (1.25, 0.5) {};
		\node [style=none] (1) at (-1.25, 0.5) {};
	\end{pgfonlayer}
	\begin{pgfonlayer}{edgelayer}
		\draw [style=swap, bend right=90, looseness=1.50] (1.center) to (0.center);
	\end{pgfonlayer}
\end{tikzpicture}

%% file: figures/boxexample3.tikz
\begin{tikzpicture}[scale=0.80]
	\begin{pgfonlayer}{nodelayer}
		\node [style=none] (0) at (-1, -1) {};
		\node [style=none] (1) at (1.5, -3) {};
		\node [style=langbox] (2) at (0.25, -0.75) {\!\!\!{\tt process 2}\!\!\!};
		\node [style=none] (3) at (1.5, -1) {};
		\node [style=none] (4) at (-1, -3) {};
		\node [style=none] (5) at (-1, -0.5) {};
		\node [style=none] (6) at (-1, -0.25) {};
		\node [style=none] (7) at (1.5, -0.5) {};
		\node [style=none] (8) at (1.5, 0.5) {};
		\node [style=langbox] (9) at (2.75, 0.75) {\!\!\!{\tt process 1}\!\!\!};
		\node [style=none] (10) at (4, 0.25) {};
		\node [style=none] (11) at (4, 0.5) {};
		\node [style=none] (12) at (-3, -0.25) {};
		\node [style=none] (13) at (-3, -3) {};
		\node [style=none] (14) at (1.5, 1) {};
		\node [style=none] (15) at (1.5, 1.75) {};
		\node [style=none] (16) at (4, 1) {};
		\node [style=none] (17) at (4, 1.75) {};
		\node [style=none] (18) at (6, 1.75) {};
		\node [style=none] (19) at (6, 0.25) {};
	\end{pgfonlayer}
	\begin{pgfonlayer}{edgelayer}
		\draw [style=swap] (1.center) to (3.center);
		\draw [style=swap] (4.center) to (0.center);
		\draw [style=swap] (5.center) to (6.center);
		\draw [style=swap] (7.center) to (8.center);
		\draw [style=swap] (10.center) to (11.center);
		\draw [style=swap] (13.center) to (12.center);
		\draw [style=swap] (14.center) to (15.center);
		\draw [style=swap] (16.center) to (17.center);
		\draw [style=swap] (19.center) to (18.center);
		\draw [style=swap, bend right=90, looseness=1.25] (10.center) to (19.center);
		\draw [style=swap, bend left=90, looseness=1.25] (12.center) to (6.center);
	\end{pgfonlayer}
\end{tikzpicture}

%% file: figures/capdiscard.tikz
\begin{tikzpicture}[scale=0.80]
	\begin{pgfonlayer}{nodelayer}
		\node [style=none] (0) at (1.75, 0) {};
		\node [style=none] (1) at (-0.75, 0) {};
		\node [style=none] (2) at (3, 0) {$=$};
		\node [style=none] (3) at (1.75, -1) {};
		\node [style=downground] (4) at (-0.75, -0.5) {};
		\node [style=upground] (5) at (4.75, 0.5) {};
		\node [style=none] (6) at (4.75, -0.75) {};
		\node [style=none] (7) at (6.75, 1.5) {};
		\node [style=none] (8) at (5.5, 1) {};
		\node [style=right label] (9) at (7, 1.5) {\color{gray} system $B$ is now maximal noise};
		\node [style=none] (10) at (-1.5, -1) {};
		\node [style=none] (11) at (-2.75, -1.5) {};
		\node [style=left label] (12) at (-3, -1.5) {\color{gray} system $A$ is being discarded};
	\end{pgfonlayer}
	\begin{pgfonlayer}{edgelayer}
		\draw [style=swap, bend left=90, looseness=1.50] (1.center) to (0.center); 
		\draw [style=swap] (0.center) to (3.center);
		\draw [style=swap] (1.center) to (4);
		\draw [style=swap] (5) to (6.center);
		\draw [style=pointer edge, in=45, out=180, looseness=1.00] (7.center) to (8.center);
		\draw [style=pointer edge, in=-135, out=0, looseness=1.00] (11.center) to (10.center);
	\end{pgfonlayer}
\end{tikzpicture}

%% file: figures/box5.tikz
\begin{tikzpicture}[scale=0.80]
	\begin{pgfonlayer}{nodelayer}
		\node [style=none] (0) at (4, -0.5) {};
		\node [style=langbox] (1) at (4, -0.25) {\!\!\!{\tt process}\!\!\!};
		\node [style=none] (2) at (4, -1.5) {};
		\node [style=none] (3) at (4, 0) {};
		\node [style=none] (4) at (4, 0.5) {};
		\node [style=none] (5) at (4, 0.5) {};
		\node [style=none] (6) at (1.5, 0.5) {};
		\node [style=none] (7) at (1.5, -1.5) {};
		\node [style=none] (8) at (0, 0) {$=$};
		\node [style=none] (9) at (-3, 0.25) {};
		\node [style=none] (10) at (-1.5, 0.25) {};
		\node [style=langstate] (11) at (-2.25, 0.5) {\!\!\!{\tt state}\!\!\!};
		\node [style=none] (12) at (-1.5, -0.75) {};
		\node [style=none] (13) at (-3, -0.75) {};
	\end{pgfonlayer}
	\begin{pgfonlayer}{edgelayer}
		\draw [style=swap] (2.center) to (0.center);
		\draw [style=swap] (3.center) to (4.center);
		\draw [style=swap, bend left=90, looseness=1.50] (6.center) to (5.center);
		\draw [style=swap] (6.center) to (7.center);
		\draw [style=swap] (12.center) to (10.center);
		\draw [style=swap] (13.center) to (9.center);
	\end{pgfonlayer}
\end{tikzpicture}

%% file: figures/box4.tikz
\begin{tikzpicture}[scale=0.80]
	\begin{pgfonlayer}{nodelayer}
		\node [style=none] (0) at (0, -0.25) {};
		\node [style=langbox] (1) at (0, 0) {\!\!\!{\tt process}\!\!\!};
		\node [style=none] (2) at (0, -1.25) {};
		\node [style=none] (3) at (0, 0.25) {};
		\node [style=none] (4) at (0, 1.25) {};
	\end{pgfonlayer}
	\begin{pgfonlayer}{edgelayer}
		\draw [style=swap] (2.center) to (0.center);
		\draw [style=swap] (3.center) to (4.center);
	\end{pgfonlayer}
\end{tikzpicture}

%% file: figures/yank.tikz
\begin{tikzpicture}[scale=0.80]
	\begin{pgfonlayer}{nodelayer}
		\node [style=none] (0) at (1.25, 0) {};
		\node [style=none] (1) at (-1.25, 0) {};
		\node [style=none] (2) at (-1.25, 0) {};
		\node [style=none] (3) at (-3.75, 0) {};
		\node [style=none] (4) at (1.25, -1.25) {};
		\node [style=none] (5) at (-3.75, 1.25) {};
		\node [style=none] (6) at (2.5, 0) {$=$};
		\node [style=none] (7) at (3.75, -1.25) {};
		\node [style=none] (8) at (3.75, 1.25) {};
	\end{pgfonlayer}
	\begin{pgfonlayer}{edgelayer}
		\draw [style=swap, bend left=90, looseness=1.50] (1.center) to (0.center);
		\draw [style=swap, bend right=90, looseness=1.50] (3.center) to (2.center);
		\draw [style=swap] (5.center) to (3.center);
		\draw [style=swap] (0.center) to (4.center);
		\draw [style=swap] (8.center) to (7.center);
	\end{pgfonlayer}
\end{tikzpicture}

%% file: figures/box3.tikz
\begin{tikzpicture}[scale=0.80]
	\begin{pgfonlayer}{nodelayer}
		\node [style=none] (0) at (0, 1.75) {};
		\node [style=langbox] (1) at (0, 3) {\!\!\!{\tt process 1}\!\!\!};
		\node [style=none] (2) at (0, 2.75) {};
		\node [style=none] (3) at (0, 3.25) {};
		\node [style=none] (4) at (0, 4.25) {};
		\node [style=langbox] (5) at (0, -2.5) {\!\!\!{\tt process $n$}\!\!\!};
		\node [style=none] (6) at (0, -2.25) {};
		\node [style=none] (7) at (0, -1.25) {};
		\node [style=none] (8) at (0, -3.75) {};
		\node [style=none] (9) at (0, -2.75) {};
		\node [style=none] (10) at (0, 1.25) {};
		\node [style=langbox] (11) at (0, 1.5) {\!\!\!{\tt process 2}\!\!\!};
		\node [style=none] (12) at (0, 0.25) {};
		\node [style=none] (13) at (0, -0.25) {$\vdots$};
	\end{pgfonlayer}
	\begin{pgfonlayer}{edgelayer}
		\draw [style=swap] (0.center) to (2.center);
		\draw [style=swap] (3.center) to (4.center);
		\draw [style=swap] (8.center) to (9.center);
		\draw [style=swap] (6.center) to (7.center);
		\draw [style=swap] (12.center) to (10.center);
	\end{pgfonlayer}
\end{tikzpicture}

%% file: figures/separate2.tikz
\begin{tikzpicture}[scale=0.80]
	\begin{pgfonlayer}{nodelayer}
		\node [style=none] (0) at (2, -2) {};
		\node [style=none] (1) at (2, 1) {};
		\node [style=none] (2) at (2, 2) {};
		\node [style=none] (3) at (-2, 1.25) {};
		\node [style=none] (4) at (-2, -0.25) {};
		\node [style=langbox] (5) at (-2, 0) {\!\!\!$f$\!\!\!};
		\node [style=none] (6) at (-2, 0.25) {};
		\node [style=langeffect] (7) at (2, 0.75) {\!\!\!\!$\psi_{A}$\!\!\!\!};
		\node [style=none] (8) at (2, -1) {};
		\node [style=langstate] (9) at (2, -0.75) {\!\!\!\!$\psi_{B}$\!\!\!\!};
		\node [style=none] (10) at (0, 0) {$\not=$};
		\node [style=none] (11) at (-2, -1.25) {};
	\end{pgfonlayer}
	\begin{pgfonlayer}{edgelayer}
		\draw [style=swap] (11.center) to (4.center);
		\draw [style=swap] (6.center) to (3.center);
		\draw [style=swap] (2.center) to (1.center);
		\draw [style=swap] (0.center) to (8.center);
	\end{pgfonlayer}
\end{tikzpicture}

%% file: figures/electricity-example2.tikz
\begin{tikzpicture}
	\begin{pgfonlayer}{nodelayer}
		\node [style=none] (0) at (-3.25, 3.25) {};
		\node [style=none] (1) at (0, -2.75) {};
		\node [style=none] (2) at (0, -4) {};
		\node [style=langbox] (3) at (0, -2.5) {\quad\qquad\ plugstrip \qquad\quad\ };
		\node [style=right label] (4) at (0.25, -3.75) {\textit{uk-plug}};
		\node [style=none] (5) at (-3.25, -2.25) {};
		\node [style=none] (6) at (3, -2.25) {};
		\node [style=none] (7) at (3, 3.25) {};
		\node [style=right label, yshift=-0.3mm] (8) at (3.25, 3) {\textit{uk-plug}};
		\node [style=none] (9) at (0, -2.25) {};
		\node [style=none] (10) at (0, -0.25) {};
		\node [style=right label, yshift=-0.3mm] (11) at (0.25, -1.25) {\textit{uk-plug}};
		\node [style=langbox] (12) at (0, 0) {adaptor};
		\node [style=none] (13) at (0, 2.25) {};
		\node [style=langstate] (14) at (0, 2.5) {power drill};
		\node [style=none] (15) at (0, 0.25) {};
		\node [style=right label, yshift=-0.3mm] (16) at (0.25, 1.25) {\textit{eu-plug}};
		\node [style=left label, yshift=-0.3mm] (17) at (-3.5, 3) {\textit{uk-plug}};
	\end{pgfonlayer}
	\begin{pgfonlayer}{edgelayer}
		\draw (2.center) to (1.center);
		\draw (0.center) to (5.center);
		\draw (7.center) to (6.center);
		\draw (10.center) to (9.center);
		\draw (15.center) to (13.center);
	\end{pgfonlayer}
\end{tikzpicture}

%% file: figures/comb.tikz
\begin{tikzpicture}
	\begin{pgfonlayer}{nodelayer}
		\node [style=langbox] (0) at (-1.5, 0) {\!\!\!\! {\tt process} \!\!\!\!};
		\node [style=none] (1) at (-1.5, 2) {};
		\node [style=none] (2) at (-1.5, -2) {};
		\node [style=none] (3) at (-1.5, 0.25) {};
		\node [style=none] (4) at (-1.5, -0.25) {};
		\node [style=none] (5) at (-1.5, -0.75) {};
		\node [style=none] (6) at (-1.5, 0.75) {};
		\node [style=none] (7) at (-2.5, 0.75) {};
		\node [style=none] (8) at (-2.5, -0.75) {};
		\node [style=none] (9) at (0.5, -0.75) {};
		\node [style=none] (10) at (0.5, 0.75) {};
		\node [style=none] (11) at (2.5, -1.25) {};
		\node [style=none] (12) at (2.5, 1.25) {};
		\node [style=none] (13) at (-2.5, 1.25) {};
		\node [style=none] (14) at (-2.5, -1.25) {};
		\node [style=none] (15) at (1.5, 0) {\footnotesize$\texttt{comb}$};
		\node [style=none] (16) at (-1.5, 1.25) {};
		\node [style=none] (17) at (-1.5, -1.25) {};
	\end{pgfonlayer}
	\begin{pgfonlayer}{edgelayer}
		\draw (4.center) to (5.center);
		\draw (6.center) to (3.center);
		\draw (7.center) to (10.center);
		\draw (8.center) to (9.center);
		\draw (10.center) to (9.center);
		\draw (13.center) to (7.center);
		\draw (13.center) to (12.center);
		\draw (8.center) to (14.center);
		\draw (14.center) to (11.center);
		\draw (12.center) to (11.center);
		\draw (1.center) to (16.center);
		\draw (17.center) to (2.center);
	\end{pgfonlayer}
\end{tikzpicture}

%% file: figures/cnot.tikz
\begin{tikzpicture}
	\begin{pgfonlayer}{nodelayer}
		\node [style=none] (0) at (-0.75, 1) {};
		\node [style=white dot] (1) at (-0.75, 0) {};
		\node [style=none] (2) at (-0.75, -1) {};
		\node [style=none] (3) at (0.75, -1) {};
		\node [style=none] (4) at (0.75, 1) {};
		\node [style=gray dot] (5) at (0.75, 0) {};
	\end{pgfonlayer}
	\begin{pgfonlayer}{edgelayer}
		\draw (1) to (0.center);
		\draw (2.center) to (1);
		\draw (5) to (4.center);
		\draw (3.center) to (5);
		\draw (1) to (5);
	\end{pgfonlayer}
\end{tikzpicture}

%% file: figures/cnotfromspiders.tikz
\begin{tikzpicture}[dotpic]
	\begin{pgfonlayer}{nodelayer}
		\node [style=none] (0) at (-1.5, -0.25) {};
		\node [style=white dot] (1) at (-0.75, -1.25) {};
		\node [style=none] (2) at (-0.75, -2.5) {};
		\node [style=none] (3) at (1.5, 0.25) {};
		\node [style=none] (4) at (0.75, 2.5) {};
		\node [style=gray dot] (5) at (0.75, 1.25) {};
		\node [style=none] (6) at (0, -0.25) {};
		\node [style=none] (7) at (0, 0.25) {}; 
		\node [style=none] (8) at (1.5, -2.5) {};
		\node [style=none] (9) at (-1.5, 2.5) {};
		\node [style=none] (10) at (-1.75, -0.25) {};
		\node [style=none] (11) at (0.25, -0.25) {};
		\node [style=none] (12) at (0.25, -2.25) {};
		\node [style=none] (13) at (-1.75, -2.25) {};
		\node [style=none] (14) at (1.75, 0.25) {};
		\node [style=none] (15) at (-0.25, 0.25) {};
		\node [style=none] (16) at (1.75, 2.25) {};
		\node [style=none] (17) at (-0.25, 2.25) {}; 
	\end{pgfonlayer}
	\begin{pgfonlayer}{edgelayer}
		\draw [in=-90, out=135, looseness=1.00] (1) to (0.center);
		\draw (2.center) to (1);
		\draw (5) to (4.center);
		\draw [in=-30, out=90, looseness=1.00] (3.center) to (5);
		\draw [style=swap, in=-90, out=45, looseness=1.00] (1) to (6.center);
		\draw [style=swap, in=-150, out=90, looseness=1.00] (7.center) to (5);
		\draw [style=swap] (7.center) to (6.center);
		\draw [style=swap] (3.center) to (8.center);
		\draw [style=swap] (9.center) to (0.center);
		\draw [style=dashed] (17.center) to (16.center);
		\draw [style=dashed] (16.center) to (14.center);
		\draw [style=dashed] (14.center) to (15.center);
		\draw [style=dashed] (17.center) to (15.center);
		\draw [style=dashed] (10.center) to (13.center);
		\draw [style=dashed] (13.center) to (12.center); 
		\draw [style=dashed] (12.center) to (11.center);
		\draw [style=dashed] (11.center) to (10.center);
	\end{pgfonlayer}
\end{tikzpicture}

%% file: figures/greenspider.tikz
\begin{tikzpicture}[scale=1.6]
	\begin{pgfonlayer}{nodelayer}
		\node [style=none] (0) at (0.75, 0.75) {};
		\node [style=none] (1) at (-0.75, -0.75) {};
		\node [style=none] (2) at (0.75, -0.75) {};
		\node [style=none, yshift=1mm] (3) at (0.25, 0.5) {...};
		\node [style=none] (4) at (-0.75, 0.75) {};
		\node [style=none] (5) at (-0.25, 0.75) {};
		\node [style=none] (6) at (-0.25, -0.75) {};
		\node [style=white phase dot] (7) at (0, 0) {$\alpha$};
		\node [style=none, yshift=-1mm] (8) at (0.25, -0.5) {...};
	\end{pgfonlayer}
	\begin{pgfonlayer}{edgelayer}
		\draw [in=-144, out=90, looseness=1.00] (1.center) to (7);
		\draw [in=-90, out=36, looseness=1.00] (7) to (0.center);
		\draw [bend left=15, looseness=1.00] (6.center) to (7);
		\draw [in=-36, out=90, looseness=1.00] (2.center) to (7);
		\draw [in=-90, out=144, looseness=1.00] (7) to (4.center);
		\draw [bend left=15, looseness=1.00] (7) to (5.center);
	\end{pgfonlayer}
\end{tikzpicture}

%% file: figures/redspider.tikz
\begin{tikzpicture}[scale=1.6]
	\begin{pgfonlayer}{nodelayer}
		\node [style=none] (0) at (0.75, 0.75) {};
		\node [style=none] (1) at (-0.75, -0.75) {};
		\node [style=none] (2) at (0.75, -0.75) {};
		\node [style=none, yshift=1mm] (3) at (0.25, 0.5) {...};
		\node [style=none] (4) at (-0.75, 0.75) {};
		\node [style=none] (5) at (-0.25, 0.75) {};
		\node [style=none] (6) at (-0.25, -0.75) {};
		\node [style=gray phase dot] (7) at (0, 0) {$\alpha$};
		\node [style=none, yshift=-1mm] (8) at (0.25, -0.5) {...};
	\end{pgfonlayer}
	\begin{pgfonlayer}{edgelayer}
		\draw [in=-144, out=90, looseness=1.00] (1.center) to (7);
		\draw [in=-90, out=36, looseness=1.00] (7) to (0.center);
		\draw [bend left=15, looseness=1.00] (6.center) to (7);
		\draw [in=-36, out=90, looseness=1.00] (2.center) to (7); 
		\draw [in=-90, out=144, looseness=1.00] (7) to (4.center);
		\draw [bend left=15, looseness=1.00] (7) to (5.center); 
	\end{pgfonlayer}
\end{tikzpicture} 

%% file: figures/flowers1rp.tikz
\begin{tikzpicture}
	\begin{pgfonlayer}{nodelayer}
		\node [style=none] (0) at (-12, 0) {};
		\node [style=none] (1) at (-8.75, 0) {};
		\node [style=none] (2) at (12, 0) {};
		\node [style=none] (3) at (8.75, 0) {};
		\node [style=langstate] (4) at (-12, 0.25) {\!\!\!{\tt Alice}\!\!\!};
		\node [style=langstate] (5) at (-4, 0.25) {\!\!\!{\tt flowers}\!\!\!};
		\node [style=langstate] (6) at (-8, 0.25) {\!\!\!{\tt likes}\!\!\!};
		\node [style=langstate] (7) at (4, 0.25) {\!\!\!{\tt Bob}\!\!\!};
		\node [style=langstate] (8) at (8, 0.25) {\!\!\!{\tt gives}\!\!\!};
		\node [style=langstate] (9) at (12, 0.25) {\!\!\!{\tt Claire}\!\!\!};
		\node [style=none] (10) at (7.25, 0) {};
		\node [style=none] (11) at (4, 0) {};
		\node [style=none] (12) at (-4, 0) {}; 
		\node [style=none] (13) at (-1.25, 0) {};
		\node [style=none] (14) at (-0.25, 0) {};
		\node [style=none] (15) at (-7.25, 0) {};
		\node [style=none] (16) at (0.75, 0) {};
		\node [style=none] (17) at (8.25, 0) {};
		\node [style=none] (18) at (-0.25, 0) {};
		\node [style=none] (19) at (-1.25, 0) {};
		\node [style=none] (20) at (-1.25, 0) {};
		\node [style=none] (21) at (-0.25, 0) {};
		\node [style=none] (22) at (-0.25, 0) {};
		\node [style=none] (23) at (-0.25, 0) {};
		\node [style=none] (24) at (0.75, 0) {};
		\node [style=none] (25) at (0.75, 0) {};
		\node [style=none] (26) at (-1.75, 0) {};
		\node [style=none] (27) at (2.25, 1) {};
		\node [style=none] (28) at (0.75, 0) {};
		\node [style=none] (29) at (0.75, 0) {};
		\node [style=none] (30) at (-1.25, 0) {};
		\node [style=none] (31) at (2.25, 0) {};
		\node [style=none] (32) at (-1.75, 1) {};
		\node [style=none] (33) at (-0.25, 0) {};
		\node [style=none] (34) at (-0.25, 0) {};
		\node [style=none] (35) at (-0.25, 0) {};
		\node [style=white dot] (36) at (-0.25, 1) {};
		\node [style=none] (37) at (0.75, 0) {};
		\node [style=none] (38) at (-0.25, 0) {};
		\node [style=none] (39) at (0.25, 1.5) {};
		\node [style=none] (40) at (-1.25, 0) {};
		\node [style=none] (41) at (0.25, 2) {\!\!\!{\footnotesize\tt that}\!\!\!};
		\node [style=none] (42) at (7.75, 0) {};
		\node [style=none] (43) at (1.75, 0) {};
		\node [style=none] (44) at (1.75, 0) {};
		\node [style=white dot, boldedge] (45) at (1.75, 0.5) {};
		\node [style=none] (46) at (1.75, 0) {};
		\node [style=none] (47) at (2, 0) {};
		\node [style=none] (48) at (1.75, 0) {};
		\node [style=none] (49) at (-8, 0) {};
		\node [style=none] (50) at (-8, -2) {};
	\end{pgfonlayer}
	\begin{pgfonlayer}{edgelayer}
		\draw [style=swap, bend right=90, looseness=1.00] (0.center) to (1.center);
		\draw [style=swap, bend right=90, looseness=1.00] (3.center) to (2.center);
		\draw [style=swap, bend right=90, looseness=1.00] (11.center) to (10.center);
		\draw [style=swap, bend right=90, looseness=1.00] (12.center) to (13.center);
		\draw [style=swap, bend right=90, looseness=0.75] (15.center) to (14.center);
		\draw [style=swap, bend right=90, looseness=0.75] (16.center) to (17.center);
		\draw [style=swap, in=180, out=90, looseness=1.00] (40.center) to (36);
		\draw [style=swap] (36) to (34.center);
		\draw [style=swap, in=90, out=0, looseness=1.00] (36) to (29.center);
		\draw [style=swap, dashed] (32.center) to (26.center);
		\draw [style=swap, dashed] (26.center) to (31.center);
		\draw [style=swap, dashed] (31.center) to (27.center);
		\draw [style=swap, dashed] (39.center) to (32.center);
		\draw [style=swap, dashed] (39.center) to (27.center);
		\draw [style=boldedge, bend right=90, looseness=0.75] (44.center) to (42.center);
		\draw [style=boldedge] (45) to (43.center);
		\draw [style=boldedge] (49.center) to (50.center);
	\end{pgfonlayer}
\end{tikzpicture}

%% file: figures/flowers1rp2.tikz
\begin{tikzpicture}
	\begin{pgfonlayer}{nodelayer}
		\node [style=none] (0) at (-0.75, 0.25) {};
		\node [style=langstate] (1) at (0, 0.5) {\!\!\!{\tt likes}\!\!\!};
		\node [style=none] (2) at (0.75, 0.25) {};
		\node [style=none] (3) at (0, 0.25) {};
		\node [style=none] (4) at (0, -0.5) {};
		\node [style=none] (5) at (0.75, -0.5) {};
		\node [style=none] (6) at (0.75, 0.25) {};
		\node [style=none] (7) at (-0.75, -0.5) {};
		\node [style=none] (8) at (-0.75, 0.25) {};
	\end{pgfonlayer}
	\begin{pgfonlayer}{edgelayer}
		\draw [style=boldedge] (3.center) to (4.center); 
		\draw [style=swap] (6.center) to (5.center);
		\draw [style=swap] (8.center) to (7.center);
	\end{pgfonlayer}
\end{tikzpicture}

%% file: figures/fake.tikz
\begin{tikzpicture}
	\begin{pgfonlayer}{nodelayer}
		\node [style=none] (0) at (-0.5, 0.25) {};
		\node [style=langstate] (1) at (0, 0.5) {\!\!\!{\tt fake}\!\!\!};
		\node [style=none] (2) at (0.5, 0.25) {};
		\node [style=none] (3) at (0.5, -0.5) {};
		\node [style=none] (4) at (0.5, 0.25) {};
		\node [style=none] (5) at (-0.5, -0.5) {};
		\node [style=none] (6) at (-0.5, 0.25) {};
	\end{pgfonlayer}
	\begin{pgfonlayer}{edgelayer}
		\draw [style=swap] (4.center) to (3.center);
		\draw [style=swap] (6.center) to (5.center);
	\end{pgfonlayer}
\end{tikzpicture}

%% file: figures/abs1.tikz
\begin{tikzpicture}[scale=0.80]
	\begin{pgfonlayer}{nodelayer}
		\node [style=none] (0) at (-4, 0.5) {};
		\node [style=langstate] (1) at (-4, 0.75) {\!\!\!{\tt Alice}\!\!\!};
		\node [style=langstate] (2) at (4, 0.75) {\!\!\!{\tt Bob}\!\!\!};
		\node [style=none] (3) at (4, 0.5) {};
		\node [style=langstate] (4) at (0, 0.75) {\ {\tt hates}\ };
		\node [style=none] (5) at (0, -1) {};
		\node [style=none] (6) at (1.5, 0.5) {};
		\node [style=none] (7) at (0, 0.5) {};
		\node [style=none] (8) at (-1.5, 0.5) {};
		\node [style=none] (9) at (4, 0.5) {};
		\node [style=none] (10) at (-1.5, 0.5) {};
		\node [style=none] (11) at (1.5, 0.5) {};
		\node [style=none] (12) at (1.5, 0.5) {};
		\node [style=none] (13) at (-1.5, 0.5) {};
		\node [style=none] (14) at (-4, 0.5) {};
		\node [style=none] (15) at (-1.5, 0.5) {};
	\end{pgfonlayer}
	\begin{pgfonlayer}{edgelayer}
		\draw [style=boldedge] (7.center) to (5.center);
		\draw [style=swap, bend right=90, looseness=1.50] (11.center) to (9.center);
		\draw [style=swap, bend left=90, looseness=1.50] (10.center) to (14.center);
	\end{pgfonlayer}
\end{tikzpicture}

%% file: figures/abs2.tikz
\begin{tikzpicture}[scale=0.52]
	\begin{pgfonlayer}{nodelayer}
		\node [style=gray dot] (0) at (1.75, 4.25) {};
		\node [style=gray dot] (1) at (1.75, -6.25) {};
		\node [style=gray dot] (2) at (1.75, 6.75) {};
		\node [style=none] (3) at (1.75, -6.25) {};
		\node [style=none] (4) at (5, 5.5) {};
		\node [style=none] (5) at (5, 4.25) {};
		\node [style=none] (6) at (5, 5.5) {};
		\node [style=none] (7) at (-4.75, -5.25) {};
		\node [style=gray dot] (8) at (1.75, -5.25) {};
		\node [style=none] (9) at (1.75, 3) {};
		\node [style=none] (10) at (1.75, 0) {};
		\node [style=none] (11) at (1.75, 4.25) {};
		\node [style=none] (12) at (5, 4.25) {};
		\node [style=white dot] (13) at (-4.75, -5.25) {};
		\node [style=none] (14) at (-4.75, 3) {};
		\node [style=none] (15) at (1.75, -5.25) {};
		\node [style=none] (16) at (-4.75, -6.75) {};
		\node [style=none] (17) at (-4.75, -5.25) {};
		\node [style=white dot] (18) at (5, 4.25) {};
		\node [style=none] (19) at (1.75, 0) {};
		\node [style=none] (20) at (5, -6.75) {};
		\node [style=none] (21) at (1.75, 3) {};
		\node [style=white phase dot small] (22) at (1.75, 3) {\tiny$\alpha$};
		\node [style=white phase dot small] (23) at (1.75, 0) {\tiny$\gamma$};
		\node [style=none] (24) at (1.75, 0) {};
		\node [style=none] (25) at (1.75, 0) {};
		\node [style=none] (26) at (1.75, 1.5) {};
		\node [style=none] (27) at (1.75, 1.5) {};
		\node [style=none] (28) at (1.75, 1.5) {};
		\node [style=none] (29) at (1.75, 1.5) {};
		\node [style=gray phase dot small] (30) at (1.75, 1.5) {\tiny$\beta$};
		\node [style=none] (31) at (1.75, -2.5) {};
		\node [style=gray phase dot small] (32) at (1.75, -2.5) {\tiny$\tilde\beta$};
		\node [style=none] (33) at (1.75, -2.5) {};
		\node [style=none] (34) at (1.75, -2.5) {};
		\node [style=none] (35) at (1.75, -4) {};
		\node [style=none] (36) at (1.75, -4) {};
		\node [style=none] (37) at (1.75, -4) {};
		\node [style=none] (38) at (1.75, -4) {};
		\node [style=none] (39) at (1.75, -2.5) {};
		\node [style=white phase dot small] (40) at (1.75, -4) {\tiny$\tilde\gamma$};
		\node [style=white dot] (41) at (1.75, -1.25) {};
		\node [style=none] (42) at (-1.5, -1.25) {};
		\node [style=none] (43) at (1.75, -1.25) {};
		\node [style=none] (44) at (1.75, -1.25) {};
		\node [style=gray dot] (45) at (-1.5, -1.25) {};
		\node [style=gray dot] (46) at (-1.5, -2.25) {};
		\node [style=none] (47) at (-1.5, -2.25) {};
		\node [style=white phase dot small] (48) at (-4.75, 3) {\tiny$\beta_A$};
		\node [style=none] (49) at (-4.75, 5.5) {};
		\node [style=gray phase dot small] (50) at (-4.75, 5.5) {\tiny$\alpha_A$};
		\node [style=gray dot] (51) at (-4.75, 6.75) {};
		\node [style=gray dot] (52) at (-1.5, 6.75) {};
		\node [style=gray phase dot small] (53) at (-1.5, 5.5) {\tiny$\alpha_{\bf p}$};
		\node [style=none] (54) at (-1.5, 5.5) {};
		\node [style=none] (55) at (-1.5, 5.5) {};
		\node [style=none] (56) at (-1.5, 5.5) {};
		\node [style=none] (57) at (-1.5, 5.5) {};
		\node [style=none] (58) at (5, 5.5) {};
		\node [style=white phase dot small] (59) at (5, 3) {\tiny$\beta_B$};
		\node [style=none] (60) at (5, 5.5) {};
		\node [style=none] (61) at (5, 5.5) {};
		\node [style=none] (62) at (5, 5.5) {};
		\node [style=none] (63) at (5, 5.5) {};
		\node [style=gray dot] (64) at (5, 6.75) {};
		\node [style=none] (65) at (5, 5.5) {};
		\node [style=none] (66) at (5, 5.5) {};
		\node [style=none] (67) at (5, 5.5) {};
		\node [style=gray phase dot small] (68) at (5, 5.5) {\tiny$\alpha_B$};
		\node [style=none] (69) at (5, 5.5) {};
	\end{pgfonlayer}
	\begin{pgfonlayer}{edgelayer} 
		\draw [style=swap] (7.center) to (14.center);
		\draw [style=swap] (17.center) to (16.center);
		\draw [style=swap, in=0, out=180, looseness=1.25] (15.center) to (7.center);
		\draw [style=swap] (8) to (1);
		\draw [style=swap] (5.center) to (4.center);
		\draw [style=swap] (12.center) to (20.center);
		\draw [style=swap, in=180, out=0, looseness=1.25] (11.center) to (5.center);
		\draw [style=swap] (0) to (2);
		\draw [style=swap] (0) to (9.center);
		\draw [style=swap] (22) to (26.center);
		\draw [style=swap] (26.center) to (25.center);
		\draw [style=swap] (31.center) to (36.center);
		\draw [style=swap, in=180, out=0, looseness=1.25] (42.center) to (44.center);
		\draw [style=swap] (25.center) to (41);
		\draw [style=swap] (35.center) to (15.center);
		\draw [style=swap] (45) to (46);
		\draw [style=swap] (51) to (49.center);
		\draw [style=swap] (64) to (60.center);
		\draw [style=swap] (52) to (53);
		\draw [style=swap] (53) to (45);
		\draw [style=swap] (49.center) to (48);
		\draw [style=swap] (41) to (39.center);
	\end{pgfonlayer}
\end{tikzpicture}

%% file: figures/flowersrpQUINT3rev.tikz
\begin{tikzpicture}
	\begin{pgfonlayer}{nodelayer}
		\node [style=none] (0) at (3.5, 2.5) {};
		\node [style=none] (1) at (-5.5, 2.5) {};
		\node [style=none] (2) at (-1.5, 2.5) {};
		\node [style=none] (3) at (-10.5, 2.5) {};
		\node [style=none] (4) at (-10.5, 3.25) {\!\!\!{\tt Alice}\!\!\!};
		\node [style=none] (5) at (-1.5, 3.25) {\!\!\!{\tt Bob}\!\!\!};
		\node [style=none] (6) at (3.5, 3.25) {\!\!\!{\tt Claire}\!\!\!};
		\node [style=none] (7) at (-5.5, 3.25) {\!\!\!{\tt flowers}\!\!\!};
		\node [style=none] (8) at (-7.25, 1.75) {};
		\node [style=white dot] (9) at (-5.5, -1) {};
		\node [style=none] (10) at (-5.5, 0.25) {};
		\node [style=none] (11) at (3.5, -2) {};
		\node [style=white dot] (12) at (-10.5, 0.25) {};
		\node [style=none] (13) at (-0.5, 0.25) {};
		\node [style=white dot] (14) at (-1.5, -1) {};
		\node [style=none] (15) at (1, 0.25) {};
		\node [style=white dot] (16) at (-5.5, 0.25) {};
		\node [style=none] (17) at (-8.75, 1.75) {};
		\node [style=none] (18) at (2.5, 0.25) {};  
		\node [style=langstate] (19) at (1, 0.5) {\!\!{\tt *gives*}\!\!};
		\node [style=langstate] (20) at (-8, 2) {\!\!\!\!{\tt *likes*}\!\!\!\!};
		\node [style=none] (21) at (1, 0.25) {};
		\node [style=none] (22) at (2.5, 0.25) {};
		\node [style=none] (23) at (-5.5, -2) {};
		\node [style=white dot] (24) at (3.5, -1) {};
		\node [style=none] (25) at (-1.5, -2) {};
		\node [style=none] (26) at (-10.5, -2) {};
	\end{pgfonlayer}
	\begin{pgfonlayer}{edgelayer}
		\draw [style=swap] (12) to (26.center);
		\draw [style=swap, in=30, out=-90, looseness=0.75] (17.center) to (12);
		\draw [style=swap, in=15, out=-90, looseness=1.00] (13.center) to (14);
		\draw [style=swap, in=45, out=-90, looseness=1.00] (22.center) to (24);
		\draw [style=swap, in=0, out=-90, looseness=1.25] (15.center) to (9);
		\draw [style=swap, in=30, out=-90, looseness=1.25] (8.center) to (10.center);
		\draw [style=swap] (3.center) to (12);
		\draw [style=swap] (1.center) to (23.center);
		\draw [style=swap] (2.center) to (25.center);
		\draw [style=swap] (0.center) to (11.center);
	\end{pgfonlayer}
\end{tikzpicture}

%% file: figures/flowersrpQUINT3revcopy.tikz
\begin{tikzpicture}
	\begin{pgfonlayer}{nodelayer}
		\node [style=none] (0) at (6.75, 5.25) {};
		\node [style=none] (1) at (-2.25, 5.25) {};
		\node [style=none] (2) at (1.75, 5.25) {};
		\node [style=none] (3) at (-7.25, 5.25) {};
		\node [style=none] (4) at (-7.25, 6) {\!\!\!{\tt Alice}\!\!\!};
		\node [style=none] (5) at (1.75, 6) {\!\!\!{\tt Bob}\!\!\!};
		\node [style=none] (6) at (6.75, 6) {\!\!\!{\tt Claire}\!\!\!};
		\node [style=none] (7) at (-2.25, 6) {\!\!\!{\tt flowers}\!\!\!};
		\node [style=none] (8) at (-4, 4.5) {};
		\node [style=white dot] (9) at (-2.25, 1.75) {};
		\node [style=none] (10) at (-2.25, 3) {};
		\node [style=none] (11) at (6.75, -2.25) {};
		\node [style=white dot] (12) at (-7.25, 3) {};
		\node [style=none] (13) at (2.75, 3) {};
		\node [style=white dot] (14) at (1.75, 1.75) {};
		\node [style=none] (15) at (4.25, 3) {};
		\node [style=white dot] (16) at (-2.25, 3) {};
		\node [style=none] (17) at (-5.5, 4.5) {};
		\node [style=none] (18) at (5.75, 3) {};
		\node [style=langstate] (19) at (4.25, 3.25) {\!\!{\tt *gives*}\!\!};
		\node [style=langstate] (20) at (-4.75, 4.75) {\!\!\!\!{\tt *likes*}\!\!\!\!};
		\node [style=none] (21) at (4.25, 3) {};
		\node [style=none] (22) at (5.75, 3) {};
		\node [style=none] (23) at (-2.25, -2.25) {};
		\node [style=white dot] (24) at (6.75, 1.75) {};
		\node [style=none] (25) at (1.75, -2.25) {};
		\node [style=none] (26) at (-7.25, -2.25) {};
		\node [style=white dot] (27) at (-7.25, 0.25) {};
		\node [style=none] (28) at (-5.5, 1.75) {};
		\node [style=white dot] (29) at (1.75, 0.25) {};
		\node [style=langstate] (30) at (-4.75, 2) {\!\!\!\!{\tt *slaps*}\!\!\!\!};
		\node [style=none] (31) at (1.75, 0.25) {};
		\node [style=none] (32) at (-4, 1.75) {};
		\node [style=none] (33) at (4.25, 0.25) {};
		\node [style=white dot] (34) at (1.75, -1.25) {};
		\node [style=langstate] (35) at (4.25, 0.5) {\!\!\!\!{\tt *cries*}\!\!\!\!};
	\end{pgfonlayer}
	\begin{pgfonlayer}{edgelayer}
		\draw [style=swap] (12) to (26.center);
		\draw [style=swap, in=30, out=-90, looseness=0.75] (17.center) to (12);
		\draw [style=swap, in=15, out=-90, looseness=1.00] (13.center) to (14);
		\draw [style=swap, in=45, out=-90, looseness=1.00] (22.center) to (24);
		\draw [style=swap, in=0, out=-90, looseness=1.25] (15.center) to (9);
		\draw [style=swap, in=30, out=-90, looseness=1.25] (8.center) to (10.center);
		\draw [style=swap] (3.center) to (12);
		\draw [style=swap] (1.center) to (23.center);
		\draw [style=swap] (2.center) to (25.center);
		\draw [style=swap] (0.center) to (11.center);
		\draw [style=swap, in=30, out=-90, looseness=0.75] (28.center) to (27);
		\draw [style=swap, in=30, out=-90, looseness=1.00] (32.center) to (31.center);
		\draw [style=swap, in=30, out=-90, looseness=0.75] (33.center) to (34);
	\end{pgfonlayer}
\end{tikzpicture}

%% file: figures/ground.tikz
\begin{tikzpicture}[scale=0.80]
	\begin{pgfonlayer}{nodelayer}
		\node [style=upground] (0) at (0, 0.5) {};
		\node [style=none] (1) at (0, -0.75) {};
	\end{pgfonlayer}
	\begin{pgfonlayer}{edgelayer}
		\draw [style=swap] (0) to (1.center);
	\end{pgfonlayer}
\end{tikzpicture}